# Внешние биллиарды

Ф.Д. Рухович
Московский физико-технический институт

**Введение**

Для любой гладкой строго выпуклой кривой на плоскости можно определить отображение внешности этой кривой в себя, называемое внешним биллиардом. А именно, обозначим кривую $\gamma$, и пусть $x$ — точка вне ее. Существуют две касательные к $\gamma$ прямые, проходящие через $x$; выберем одну из них, например правую относительно $x$, и, отразив $x$ относительно точки касания, получим новую точку $Tx$:

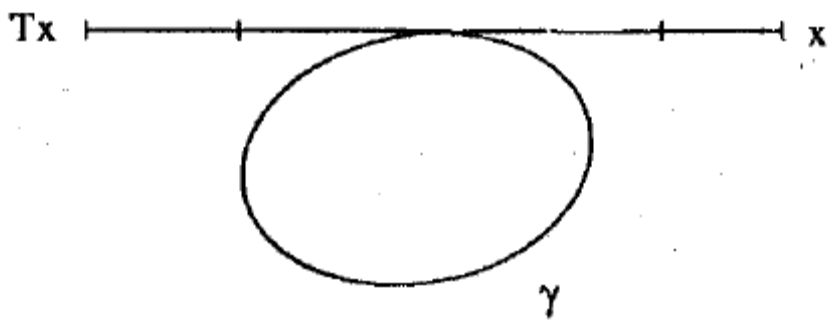

Рис.1 Определение внешнего биллиарда

Отображение $T$ называется *внешним биллиардом;* кривая $\gamma$ называется кривой *внешнего биллиарда.* Множество точек $T^n x$, где n – неотрицательное целое число, будем называть траекторией, или орбитой, внешнего биллиарда точки $x$ относительно стола $\gamma$.

В данной работе столами являются многоугольники, преимущественно правильные. В этом случае, точка правого касания определена неоднозначно, если она лежит на соответствующем продолжении одной из сторон стола; будем считать, что для таких точек преобразование внешнего биллиарда не определено.

Точку x вне фигуры $\gamma$ назовем периодической, если существует такое натуральное n, что $T^n x = x$, а периодом этой точки – минимальное такое n. В случае многоугольных столов, точки можно разбить на три типа: конечные (точки с конечной траекторией, т. е. такие точки x, что для некоторого целого неотрицательного n $T^n x$ не определено), периодические и апериодические (точки с бесконечной апериодической траекторией; такие точки мы иногда будем называть бесконечными).

**Целью** данной работы являются исследование внешних биллиардов вокруг правильных многоугольников.

**Задачами** работы являются обнаружение самоподобных структур, а также поиск точек с бесконечной апериодической траекторией для внешних биллиардов вокруг правильных восьми- и двенадцатиугольника.

Внешние биллиарды интересны своей наглядностью. В данной работе почти не будет формул, а значительной частью доказательств будут являться картинки. Такая наглядность действительно редко встречается в современной математике, как видится автору. Кто-то может заявить об «игрушечности» этой темы – да, это всего лишь игра. Но разве математика не является одной большой прекрасной игрой? И есть ли лучшее средство для развития, чем игра? К тому же одним из наиболее общих методов исследования в науке является метод построения игрушечных моделей с последующим их усовершенствованием этих моделей. За примером далеко ходить не надо: сам внешний биллиард на заре своего существования рассматривался Мозером как игрушечная модель движения планет, ибо орбита точки вокруг стола внешнего биллиарда напоминает орбиту небесного тела. Как и в случае планетарных движений, динамику двойственного биллиарда легко определить, но трудно проанализировать; в частности, совсем не ясно, может ли орбита точки уйти на бесконечность или же «упасть» на стол(цит. по [2]); этот вопрос был первоначально поставлен Б.Нейманом, который, по-видимому, и ввёл внешние (или «двойственные») биллиарды в конце 1950-х годов.

В той же книге [2] Табачников приводит две мотивации к изучению внешних биллиардов. Приведем их, почти без изменений, и мы.

«Начнем с двух мотиваций. Сначала… дадим интерпретацию двойственной биллиардной системы как механической системы, а именно импульсного осциллятора… Рассмотрим гармонический осциллятор на прямой, то есть частицу, координата которой, как функция времени, есть линейная комбинация $\sin t$ и $\cos t$. Имеется также $2\pi$-периодически движущаяся массивная стена слева от частицы, положение которой $p(t)$ удовлетворяет дифференциальному уравнению $p''(t) + p(t) = r(t)$, где $r(t)$ — это неотрицательная периодическая функция, которая удовлетворяет условиям

$$\int_0^{2\pi} r(t) \sin t \, dt = \int_0^{2\pi} r(t) \cos t \, dt = 0.$$

Когда частица ударяется о стену, происходит упругое отражение так, что скорость частицы относительно стены мгновенно меняет знак.

Эта механическая система изоморфна двойственному биллиарду относительно замкнутой выпуклой кривой $\gamma(t)$, параметризованной углом, образованным касательной с горизонтальным направлением, радиус кривизны которой г(t). Выберем начало координат $O$ внутри $\gamma$ и пусть $p(t)$ — опорная функция. Как мы уже знаем… $p''(t) + p(t) = r(t)$. Пусть $x$ — точка вне $\gamma$ и пусть плоскость вращается с постоянной угловой скоростью относительно начала координат $O$. Рассмотрим проекции $x$ и $\gamma$ на горизонтальную прямую. Положение вращающейся точки определяется как функция времени $t$ соотношением ($R\cos(t+t_0)$, $R\sin(t+t_0)$). Следовательно, проекция этой точки $x$ есть гармонический осциллятор на прямой, правая концевая точка проекции $\gamma$ есть «стена» $p(t)$. Когда осциллятор и стена соударяются, касательная из $x$ в $\gamma$ будет вертикальной. Для того чтобы в проекции было упругое отражение, точка $x$ должна отражаться от точки касания (см. рис. 2).

Второй мотивацией и объяснением термина «двойственный биллиард» является сферическая двойственность... Напомним, что на единичной сфере имеет место двойственность между точками и ориентированными прямыми (то есть большими окружностями): полюсу соответствует ориентированный экватор (см. рис. 3). Заметим, что сферическое расстояние $AB$ равно углу между линиями a и b.

Как и проективная двойственность…, сферическая двойственность распространяется на гладкие. кривые: кривая $\gamma$ определяет однопараметрическое семейство касательных, а каждая прямая определяет двойственную точку. Результирующее однопараметрическое семейство точек образует двойственную кривую $\gamma^*$…

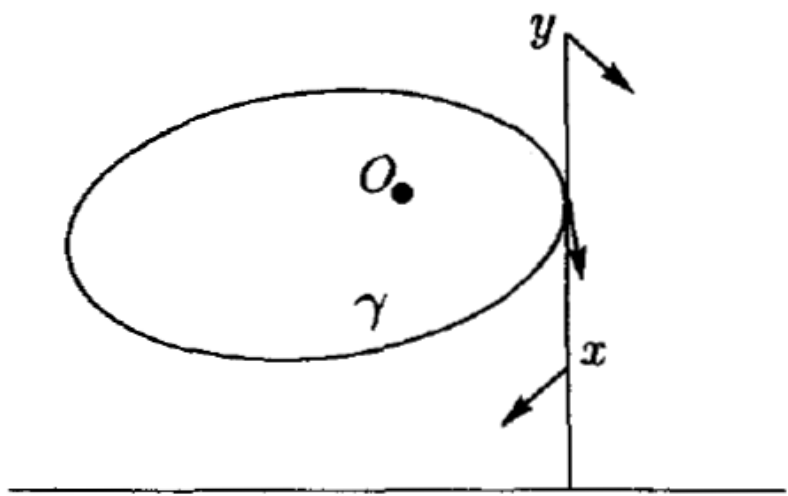

Рис. 2. Двойственный биллиард как импульсный осциллятор

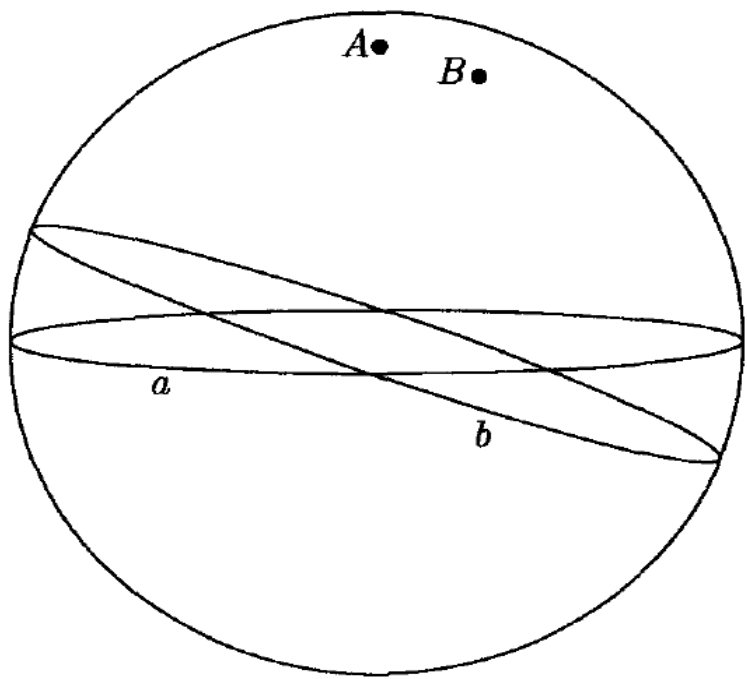

Рис.3. Сферическая двойственность

Рассмотрим биллиардное отражение от кривой *γ* (см. рис. 4). Закон биллиардного отражения читается: угол падения равен углу отражения. В терминах двойственной картины это означает, что *AL = LB,* и, следовательно, двойственное биллиардное отражение относительно двойственной кривой *γ \** переводит *A* в *B*. Таким образом, внутренний и внешние биллиарды сопряжены сферической двойственностью и две системы изоморфны на сфере. В плоскости внутренний и внешний биллиарды не зависят друг от друга так непосредственно и не существует прямой связи между системами … »

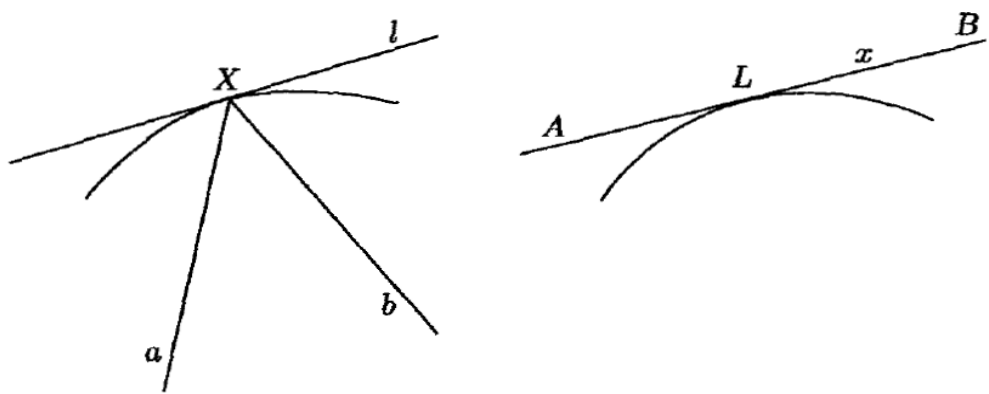

Рис. 4. Двойственность между внутренним и внешним биллиардами.

## Глава 1. Литературный обзор

Как уже было замечено во введении, внешние биллиарды стали популярны благодаря Ю.Мозеру. Он же ([4, 5]) письменно зафиксировал поставленный устно Нейманом вопрос: может ли траектория внешнего биллиарда уходить на бесконечность, т. е. быть неограниченной?

Первым, неформальным в строгом математическом смысле, ответом оказались компьютерные эксперименты, поставленные в начале девяностых годов XX века, для случая, когда столом является полуокружность. На Рис.5а черным цветом выделены точки разрыва отображения T; эти точки «окружают» большие области, содержащие точки с нечетным периодом; помимо этого, на рисунке видны и более мелкие области, которые, согласно экспериментальным данным, состоят из уходящих по спиралевидной траектории на бесконечность. Эти данные позволили Табачникову и Монро [6] предположить, что в случае полуокружности существует неограниченная траектория; одна из таких траекторий изображена на рис. 5б. Однако строгое в математическом смысле доказательство существования неограниченной траектории внешнего биллиарда относительно полуокружности удалось получить лишь в 2009 году [7].

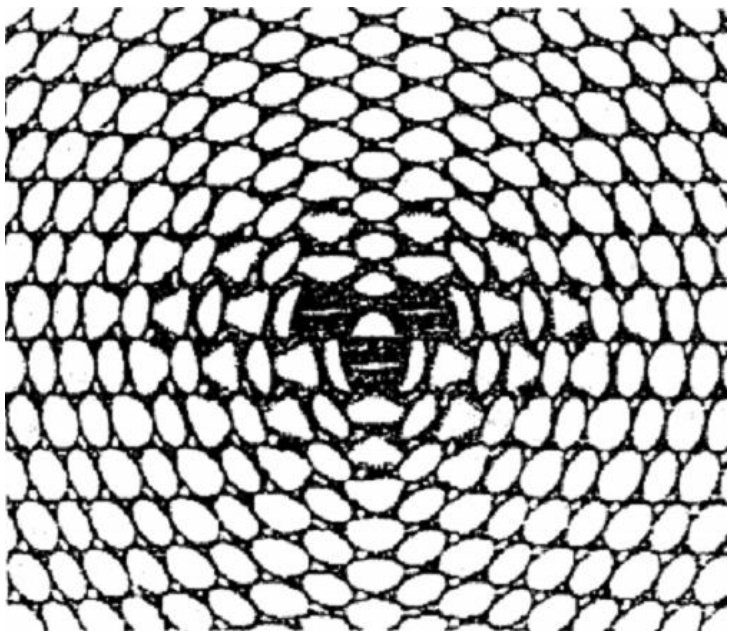

Рис. 5а. Внешний биллиард вокруг полуокружности: точки разрыва

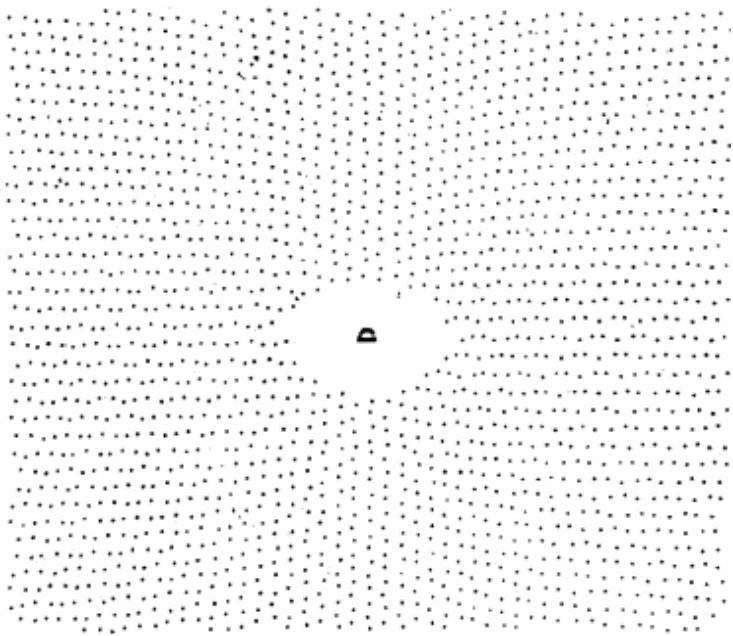

Рис. 5б. Внешний биллиард вокруг полуокружности: неограниченная траектория

Неограниченная траектория может существовать и для внешних биллиардов относительно многоугольников. Так, в 2007 году Р.Шварц доказал, что таким свойством обладает внешний биллиард вокруг «воздушного змея»(kite) — четырехугольника, вершины которого имеют координаты (0, -1), (-1, 0), (0, 1), (0, $x$), где $x$ — некоторое иррациональное число, принадлежащее интервалу (0, 1). Заметим, что ограничение на величину $x$ здесь несущественно, ибо преобразование внешнего биллиарда коммутирует с (сохраняющими ориентацию) аффинными преобразованиями; следовательно, аффинное преобразование сохраняет такие свойства (или антисвойства) орбиты, как конечность, периодичность и ограниченность; так как воздушные змеи с различными иррациональными $x$ аффинно-эквивалентны, то ограничение $x < 1$ можно снять.

Таким образом, утверждение об ограниченности траекторий внешнего биллиарда в общем случае неверно; тем не менее, найдены важные частные случаи внешних биллиардов, в которых все орбиты ограничены. Чтобы описать эти частные случаи внешних биллиардов, поговорим немного в динамике внешних биллиардов на бесконечности. Пусть точка $x$ далека от стола $\gamma$, тогда прямые $(x, Tx)$ и $(Tx, T^2x)$ практически параллельны; зафиксируем точку $x$, а кривую $\gamma$ будем сжимать посредством гомотетии к некоторой точке О внутри стола $\gamma$. Тогда сдвиг $(x, T^2x)$ будет стремиться к нулю, а траектория будет стремиться к некоторой непрерывной замкнутой кривой Г(рис. 6).

Например, в случае, когда стол является квадратом, Г также становится квадратом, повернутым на 45 градусов относительно стола; в случае треугольника и шестиугольника Г есть шестиугольник, стороны которого параллельны сторонам фигуры. Для полукруга же предельной кривой являются две параболы, пересекающиеся под прямым углом (рис.7).

Сейчас мы попробуем понять, как должна выглядеть предельная траектория Г для многоугольника. Рассмотрим произвольный выпуклый многоугольник. Плоскость вне его разбивается на «части», т.ч. $T^2$ есть композиция симметричных отражений с центрами в двух вершинах, фиксированных для каждого из множеств. Т.к. описанная композиция является параллельным переносов, то при $T^2$ каждая из «частей» будет сдвигаться на фиксированный вектор, равный по длине удвоенной диагонали (либо стороне) многоугольника.

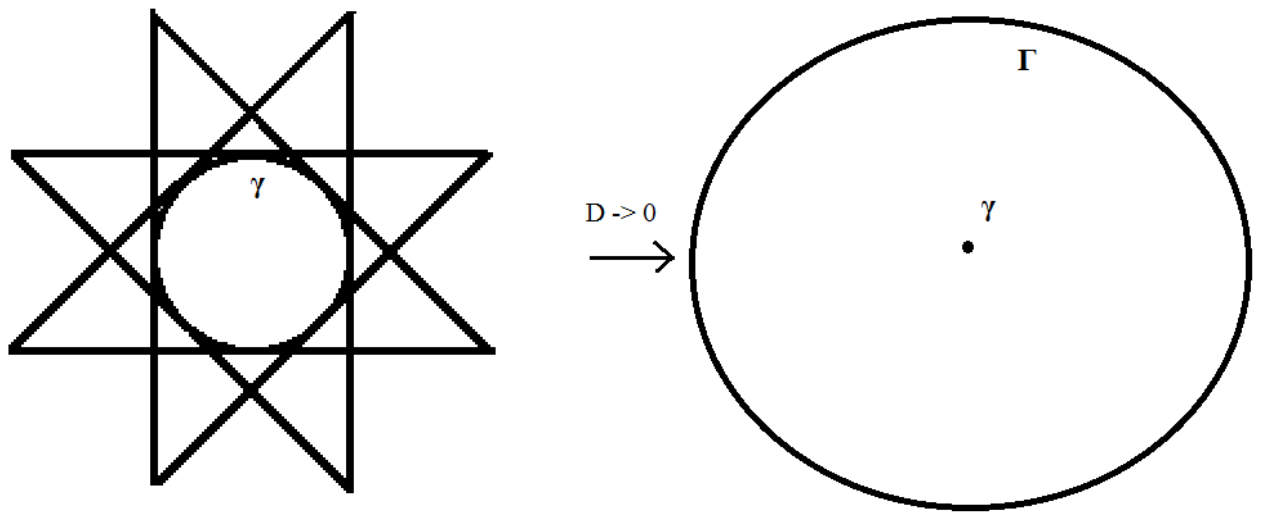

Рис. 6. Предельный переход для случая окружности

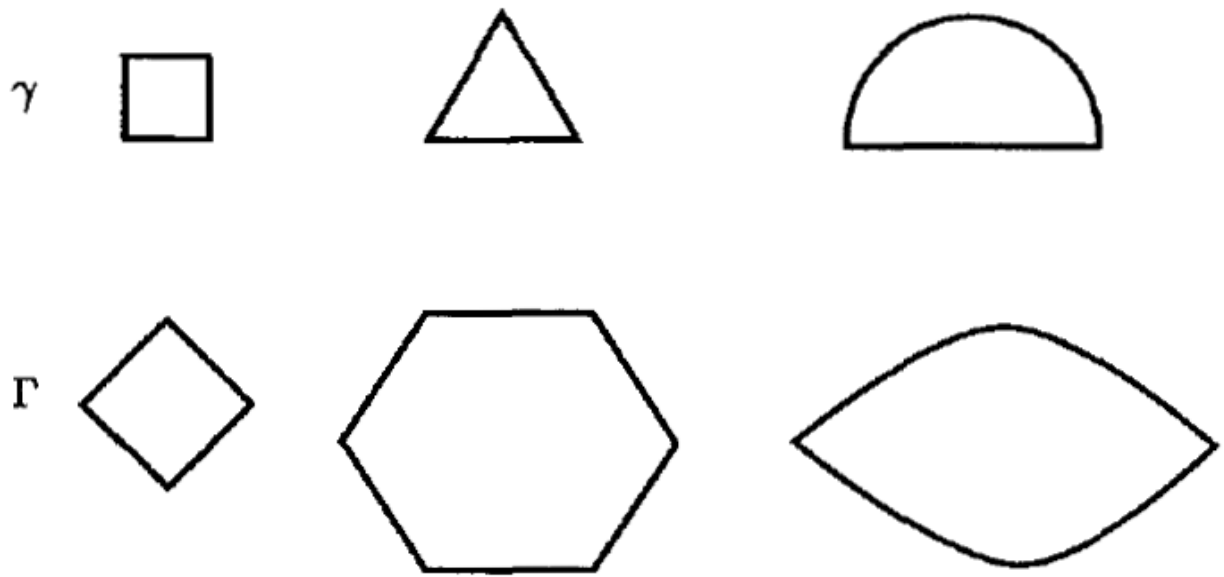

Рис. 7. Предельные кривые (экспериментальные данные)

В качестве примера рассмотрим преобразование $T^2$ вокруг квадрата,

изображенное на рис.8. При уменьшении многоугольника конечные «части» и «полосы», будут «схлопываться», в пределе исчезая. В остальных же «частях» направление сдвига меняться не будет, и в пределе Г будет двигаться по каждой такой «части» параллельно соответствующей диагонали, что и приводит нас к имеющемуся результату. В случае квадрата Г будет замкнутой. Верно ли это в случае произвольного многоугольника?

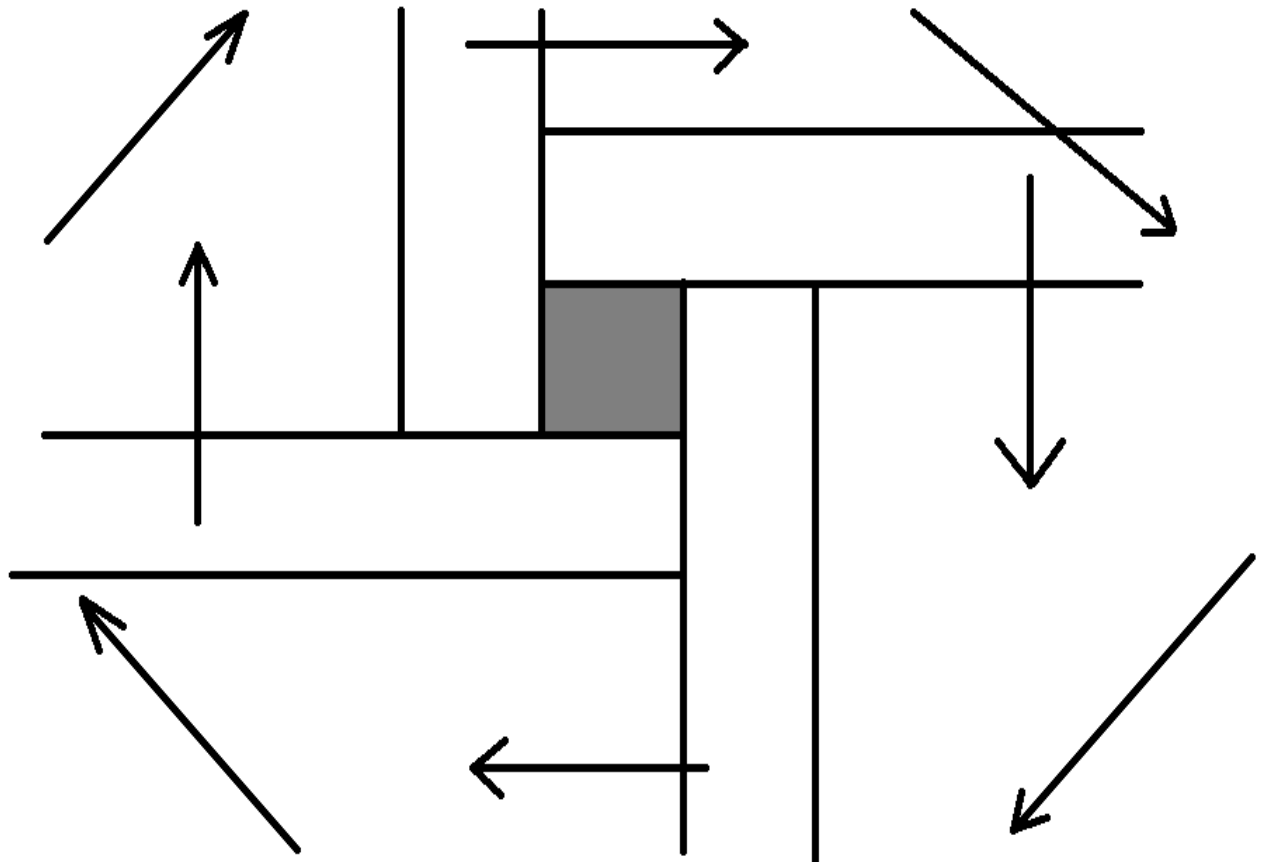

Рис.8. $T^2$ для квадрата

Чтобы ответить на этот вопрос, попробуем найти какое-нибудь выражение Г. Будем рассуждать в «полуфизическом» стиле. Пусть γ – строго выпуклая гладкая кривая, и точка A = A(t) лежит на Г. Рассмотрим две касательные к γ, параллельные OA(t). Заметим, что угол между этими касательными и прямыми (A, T(A)) и (T(A), $T^2$(A)) стремится к нулю. Следовательно, скорость точки должна быть пропорциональна *v(t)*, где *v(t)* – вектор, соединяющий точки касания параллельных касательных (рис. 9).

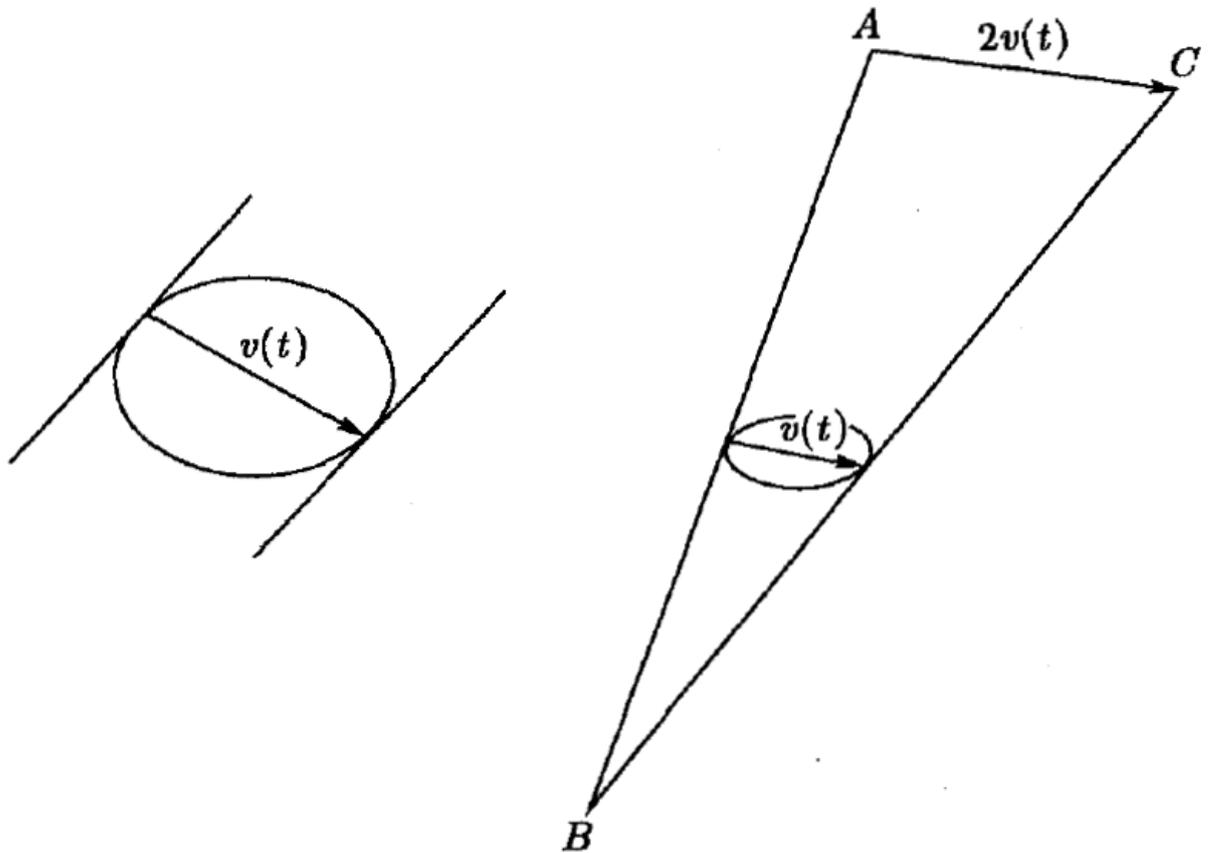

Рис.9. Определение *v(t)*

Итак, Γ задается уравнением вида

$$\Gamma'(t) = a(t)v(t),$$

где *a(t)* – некая числовая функция. Т.к. в каждой точке направление задано, то решение существует и единственно (с точностью до гомотетии). Этим решением является, как несложно проверить,

$$\Gamma(t) = \frac{v'(t)}{[v(t), v'(t)]},$$

где [..., ...] – векторное произведение. Заметим, что точка А движется со скоростью 2v(t), а *[v(t), Γ(t)] = 1*. Следовательно, выполняется второй закон Кеплера: вектор ОГ(t) заметает одинаковую площадь в единицу времени. Этот же закон верен и для многоугольника. Следовательно, Г замкнута, ибо наличие двух различных сонаправленных векторов ОГ($t_1$) и ОГ($t_2$) и равенство v($t_1$) и v($t_2$) противоречат второму закону Кеплера. Также очевидна и центральная симметричность кривой Г.

Итак, нам удалось в первом приближении описать асимптотическое поведение второй итерации внешнего биллиардного отображения на бесконечности. Разумеется, это лишь приближение, но оно имеет важные следствия [далее цит. по [2] ].

Первое следствие заключается в том, что если кривая γ достаточно гладкая (по [2], достаточно $C^5$), то, используя теорию Колмогорова – Арнольда – Мозера, можно показать, что преобразование внешнего биллиарда вне γ имеет сколь угодно далекие от биллиардного стола инвариантные кривые, которые обеспечивают ограниченность биллиардных траекторий.

Второе следствие касается случая, когда биллиардный стол – выпуклый многоугольник. В этом случае, как уже было выяснено, Г – центрально-симметричный выпуклый многоугольник (определенный с точностью до гомотетии. Каждой стороне Г соответствует «время» $t_i$, за которое точка «проходит» эту сторону; набор чисел t = ($t_1$, $t_2$, …, $t_k$) определен с точностью до общего множителя. Внешний биллиард называется квазирациональным, если t можно сделать целочисленным вектором. Доказана (см. [11, 12, 13])

**Теорема 1**: *Если внешний биллиардный стол – квазирациональный многоугольник, то ни одна орбита внешнего биллиардного соображения не уходит на бесконечность.*

Также заметим, что квазирациональными многоугольниками являются, в частности, решеточные многоугольники (т.е. многоугольники, аффинно-эквивалентные многоугольникам с рациональными координатами). Таким образом, траектории внешних биллиардов относительно решеточных и правильных многоугольников ограничены; в частности, ограничены все траектории и внешних биллиардов относительно "воздушных змеев" Шварца с рациональными $x$.

Стоит сказать несколько слов о проблеме наличия *периодической* траектории вокруг различных столов. Разберем этот вопрос подробно для случая строго выпуклого стола γ. В этом случае верна

**Теорема 2**: для любой гладкой замкнутой выпуклой кривой γ, т.ч. не существует (невырожденных) отрезков, лежащих целиком на кривой γ, и для любого натурального числа n ≥ 3 существует периодическая точка с периодом n.

**Доказательство** проведем в «полуфизическом» стиле. Рассмотрим n-угольник минимальной площади, описанный вокруг фигуры γ. Пусть это многоугольник $A_1A_2…A_n$, а отрезок $A_1A_2$ касается γ в точке B. Пусть B' – бесконечно близкая к B точка кривой, находящаяся «ближе» к точке $A_2$. Проведем через B' касательную к γ, и пусть эта касательная пересекает лучи $A_nA_1$ и $A_3A_2$ в точках $A_1'$ и $A_2'$ соответственно (рис.10).

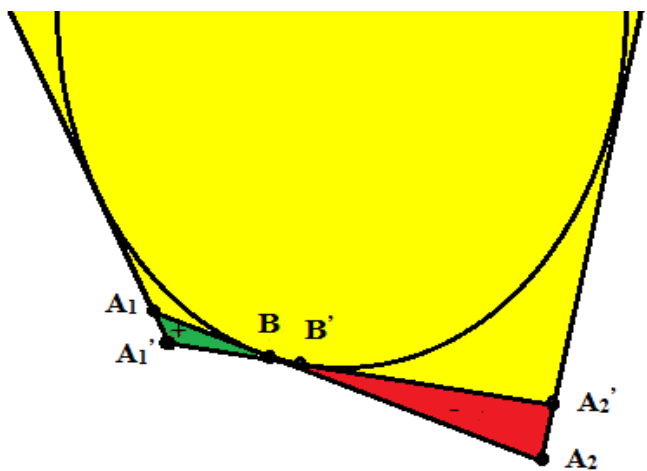

Рис.10. Изменение касательной на бесконечно малую величину

Тогда в первом порядке малости имеем

$\Delta S := S_{A_1'A_2'A_3...A_n} - S_{A_1A_2...A_n} \approx \frac{1}{2} * d\alpha * (|A_2B|^2 - |A_1B|^2)$, где $d\alpha$ – угол между прямыми $A_1A_2$ и $A_1'A_2'$. В силу минимальности площади $A_1A_2...A_n$ имеем $\Delta S \geq 0$, откуда $|A_2B| \geq |A_1B|$. Аналогично можно показать, что $|A_1B| \geq |A_2B|$; следовательно, $|A_2B| = |A_1B|$. Таким образом, в минимальном по площади описанном многоугольнике точки касания делят стороны пополам; следовательно, вершины этого многоугольника суть траектория внешнего бильярда, что и требовалось доказать.

Более того, аналогичным способом можно показать, что существуют и траектории нужного периода, «обходящие» фигуру не один раз, а любое нужное нам число раз! Соответствующая теорема выглядит следующим образом:

**Теорема 3**: для любой гладкой замкнутой выпуклой кривой γ, т.ч. не существует (невырожденных) отрезков, лежащих целиком на кривой γ, и для любых натуральных чисел n и k, т.ч. n >= 2*k+1 существует периодическая точка с периодом n, причем её траектория «обходит» фигуру ровно k раз.

Теоремы 2 и 3 неверны в случае, если стол γ есть выпуклый многоугольник. Например, как мы увидим позже, если стол γ является квадратом, то периодами для γ являются лишь числа, делящиеся на 4. Однако, из теорем следует, что внешние биллиарды вне строго выпуклых фигур обладают периодическими траекториями. Обладают ли таким свойством многоугольные столы?

Для внешнего биллиарда вокруг правильного n-угольника орбита строится тривиальным образом. В случае же решетчатого многоугольника заметим, что не существует бесконечных апериодических траекторий. Действительно, пусть, без ограничения общности, стол есть многоугольник с целыми координатами, и пусть точка *P* с координатами *(x, y)* лежит вне стола, и ее траектория не является конечной. Тогда преобразование $T^2$ сохраняет дробные части *x* и *y,* т.е. «четная часть» орбиты P находится на некотором «сдвиге» целочисленной решетки; так как орбита P ограничена по теореме 1, то траектория P периодична, что и требовалось показать. В силу же того, что размерность Хаусдорфа всех точек с конечной траекторией есть 1, то для решеточного многоугольника почти все точки являются точками с периодической траекторией. Существует доказательство и более общего факта [14]: для любого выпуклого многоугольника существует периодическая траектория. Вопрос о существовании же периодической траектории внешнего биллиарда для произвольного выпуклого стола, по-видимому, остается открытым.

Еще одной открытой проблемой является основная для нашего исследования проблема существования бесконечной апериодической траектории вокруг правильного n-угольника. В случаях n = 3, 4, 6 стол является решеточным; следовательно, согласно доказанному ранее, бесконечных апериодических траекторий не существует. В случае же n = 5 в 1993 г. Табачников [1] провёл детальнейший анализ «первой компоненты» и показал, что размерность Хаусдорфа множества непериодических точек (т. е. точек с конечной или бесконечной апериодической траекториями) есть $\log(6) / \log(\sqrt{5} + 2)$ - число, большее единицы; следовательно, точки с бесконечной апериодической траекторией существуют и обладают той же размерностью (рис.10).

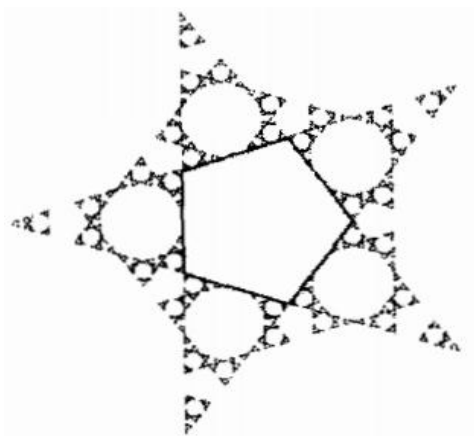

Рис.10: одна из апериодических траекторий вокруг пятиугольника

В статье [3] Табачников приводит результаты компьютерных экспериментов, дающие возможность предположить существование точки с бесконечной апериодической траекторией для правильных n-угольников с n > 6 (рис. 11-13), однако строгого анализа до сих пор нет. Таким образом, вопрос о существовании точки с бесконечной апериодической траектории для внешнего биллиарда относительно правильного n-угольника при n > 6 остается открытым.

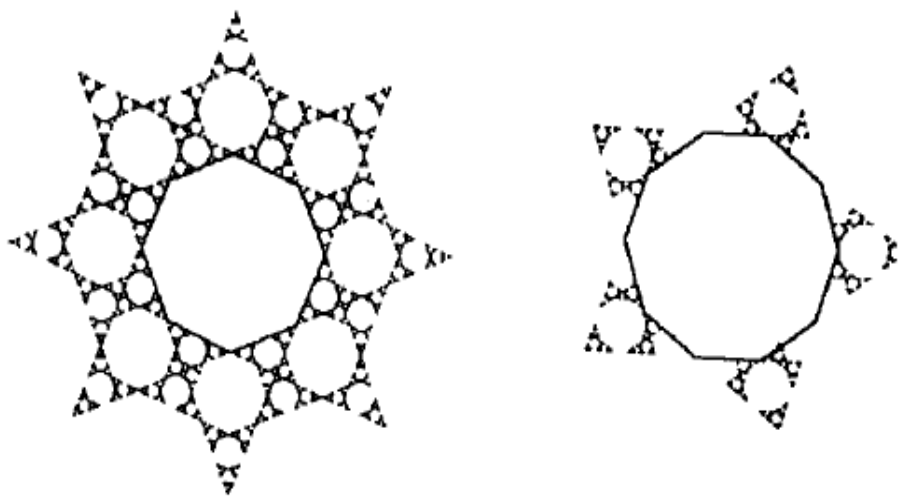

Рис.11. Эксперименты для восьми- и десятиугольника

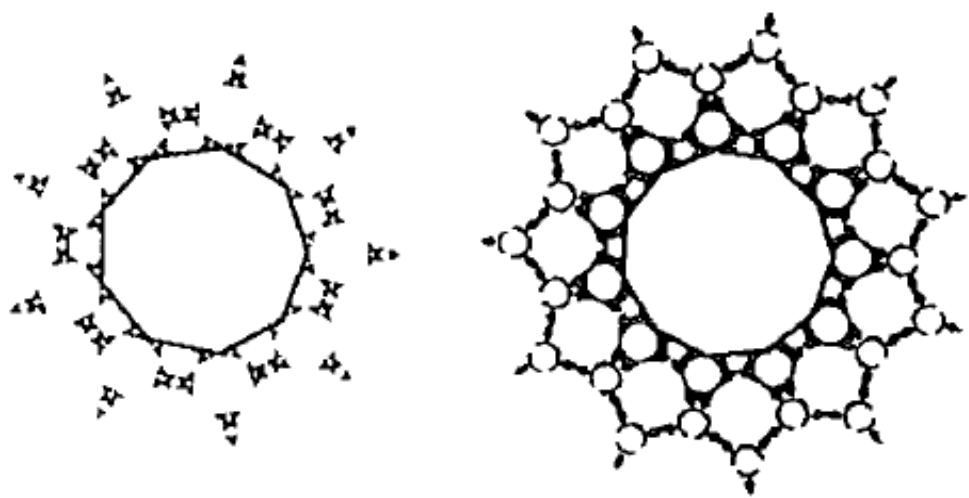

Рис.12. Эксперименты для девяти- и одиннадцатиугольника

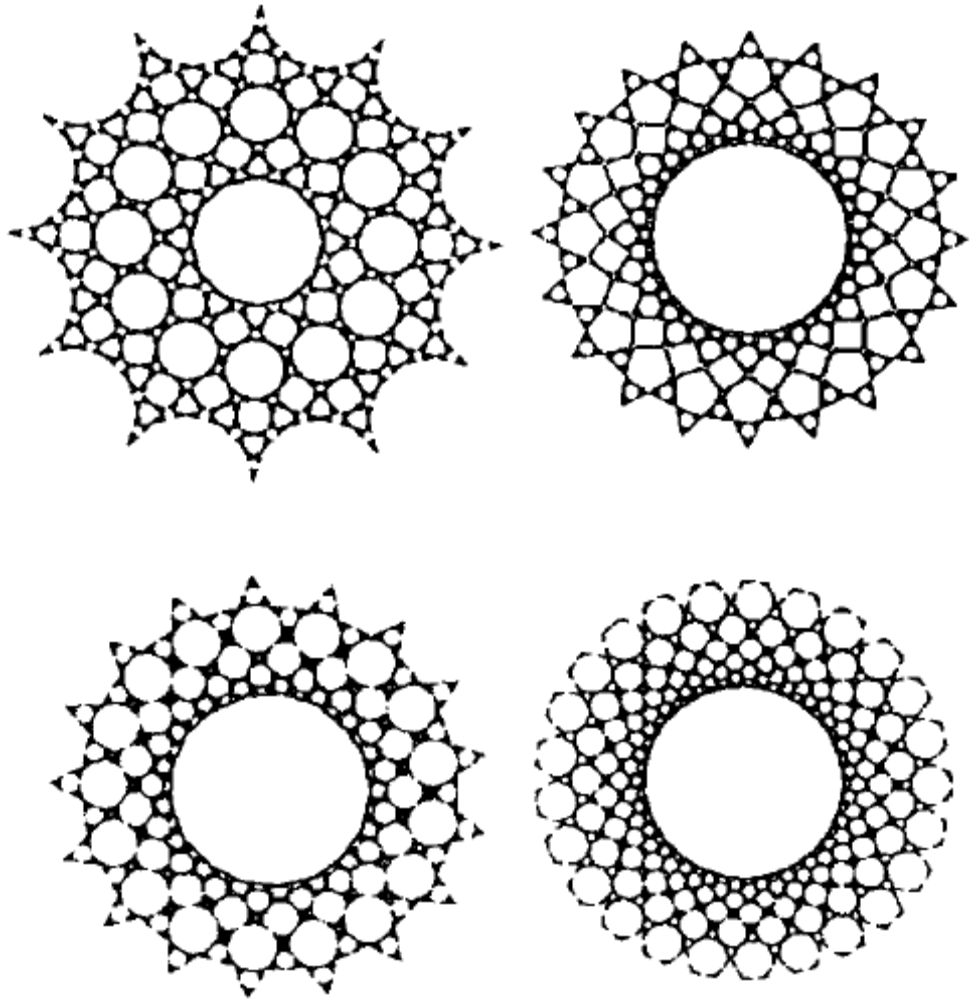

Рис.13 Эксперименты для фигур с большим числом углов

Следующим по сложности исследования является, по-видимому, правильный восьмиугольник. В своей большой статье [10] Р. Шварц проводит большое исследование внешнего биллиарда относительно восьмиугольника; в частности, он провел анализ преобразования, задающегося следующим образом. Рассмотрим правильный восьмиугольник и многоугольник X, заданный на рис. 14. Этот многоугольник можно разбить двумя способами на два равнобедренных прямоугольных треугольника $A_1$ и $B_1$ либо $A_2$ и $B_2$, изображенные на рис.15. Существует и единственно движение $\varphi_A$, переводящее $A_1$ в $A_2$; это преобразование является поворотом на $3\pi/4$ радиан вокруг центра маленького восьмиугольника на рис.14. Аналогично, существует единственное движение $\varphi_B$, переводящее $B_1$ в $B_2$, являющееся поворотом на π/4 радиан вокруг центра большого восьмиугольника. Преобразование φ: X → X Шварц определяет как объединение ограничений $φ_A$ на $A_1$ и $φ_B$ на $B_1$; такое преобразование определено корректно всюду, кроме общего отрезка треугольников $φ_A$ и $φ_B$.

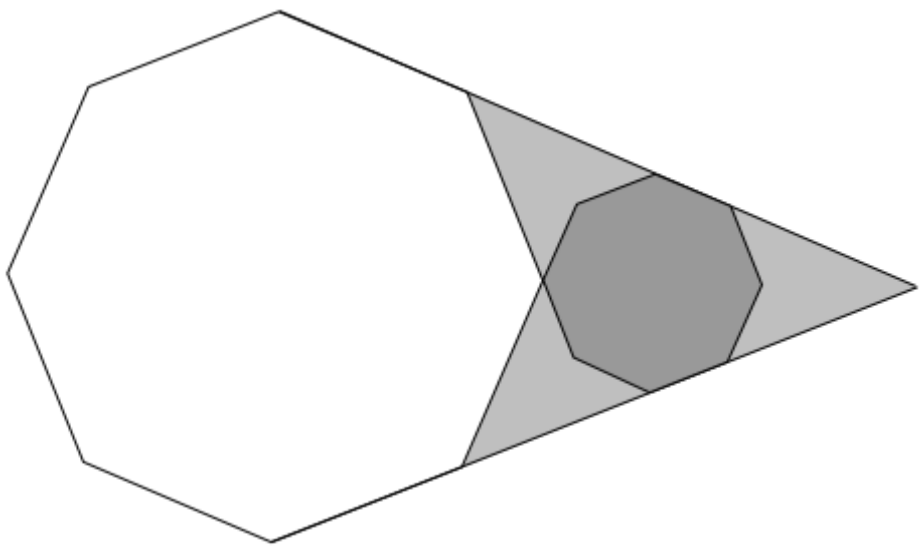

Рис.14. Многоугольник X

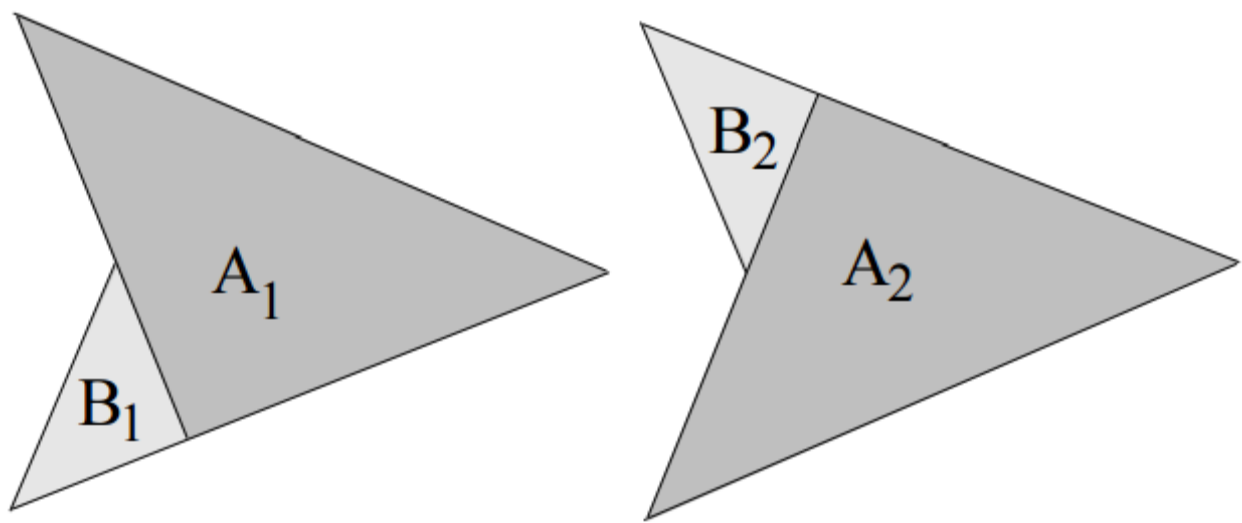

Рис.15. Два разбиения фигуры X

Именно для такого преобразования Шварцу удается получить результаты, аналогичные полученным Табачниковым для правильного восьмиугольника; в частности, ему удается получить т.н.

«ренормализационную схему»(The Renormalization Scheme) и соответствующее самоподобное преобразование, с помощью которых удается классифицировать все периодические компоненты (многоугольники с одинаковой «орбитой») и показать существование бесконечной апериодической траектории. Однако преобразование φ и четырехугольник X не связаны напрямую с внешним биллиардом вокруг восьмиугольника и рассматриваются Шварцем исключительно в качестве игрушечной модели (Toy Model). Таким образом, вопрос о существовании бесконечной апериодической траектории для внешнего биллиарда вокруг правильного восьмиугольника остается открытым.

Хотя напрямую связанная с внешним биллиардом «ренормализационная схема» описана лишь в случае правильного пятиугольника ([1]), авторы предполагают, что такая схема существует. Так, в той же статье [10] Шварц рассуждает о том, что внешние биллиарды для правильных n-угольников в случаях n = 5, 10, 8, 12, очень похожи друг на друга; в частности, в этих случаях Шварц предполагает существование эффективных ренормализационных схем, которые дают возможность получить («как минимум в принципе») полное представление о поведении преобразования внешних биллиардов. Выбор именно таких четырех чисел обусловлен особыми свойствами этих чисел; в частности, функция Эйлера, а вместе с ней и степень расширения поля $[Q(\exp(2\pi i/n)) : Q]$ равна 4 для этих, и только для этих, чисел. Заметим, что случаи n = 8 и 12 также оказываются наиболее удобными (после 3, 4 и 6) случаями и для точных компьютерных исследований, ибо координаты всех возникающих в этих случаях точек, включая вершины многоугольника, координаты первой компоненты и т. п. лежат в полях $Q(\sqrt{2})$ и $Q(\sqrt{3})$ соответственно. Однако повторимся, что существование эффективных ренормализационных схем для случаев n > 6, до сих пор не доказано, равно как и существование апериодических траекторий внешних биллиардов в этих случаях.

Что касается остальных случаев, то компьютерные исследования Шварца, изложенные в [15], обнаруживают существенно более сложное поведение динамики внешнего биллиарда. Например, в случае n = 7 Шварц обнаружил совершенно нетипичные для случаев n = 5, 8, 10, 12 периодические компоненты; например, на рис. 16 изображены периодические пятиугольники, имеющие период 57848 и диаметр порядка 0.0003 (при радиусе семиугольника 1); также Шварцу обнаружить и гораздо более «странные» периодические фигуры, обладающими периодами до 1048756, изображенные на рис.17. Из своих экспериментальных данных Шварц делает вывод о том, что полной ренормализационной схемы для внешнего биллиарда относительно правильного семиугольника не существует; также Шварц полагает, что вопрос о существовании бесконечного числа неподобных периодических компонент напрямую связан с неразрешенным на текущий момент вопросом на тему того, является ли ограниченным множество целых чисел, появляющихся в записи числа $e^{2\pi i/7}$ в цепную дробь.

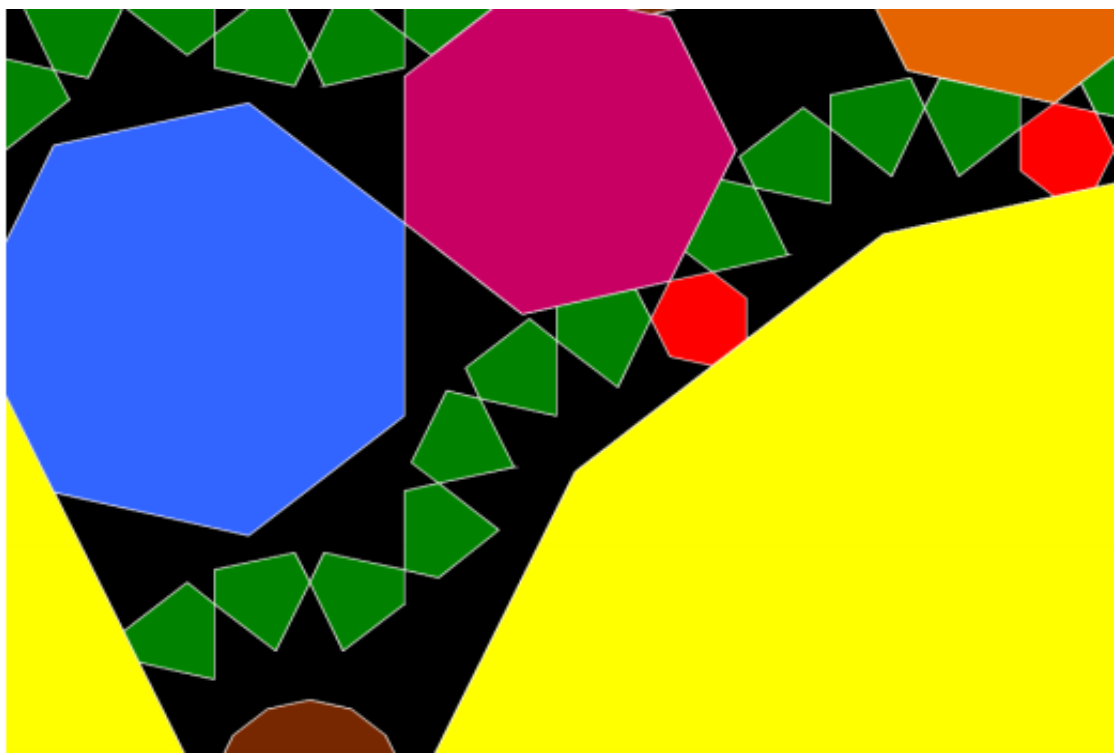

Рис.16. Асимметричные периодические компоненты внешнего биллиарда вокруг правильного семиугольника

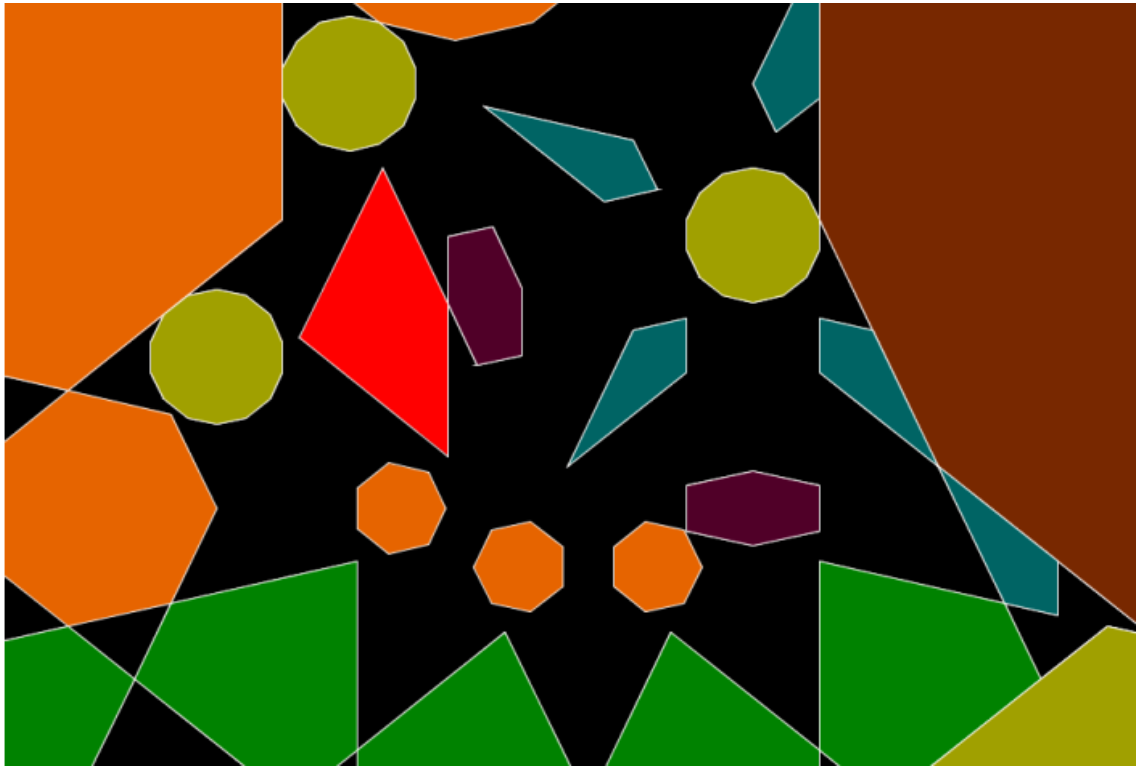

Рис.17. Некоторые «экзотические» периодические компоненты внешнего биллиарда вокруг правильного семиугольника

Исследования вокруг внешнего биллиарда ведутся и в других направлениях. Например, в (уже упомянутой) недавней работе [15] Ин-Джи Ёнг(In-Jee Jeong) проводит исследования внешнего биллиарда со сжатием (outer billiard with contraction), т. е. преобразования $T_\lambda$, $0 < \lambda <= 1$, строящегося следующим образом. Пусть дан стол $\gamma$ и точка $x$ вне стола. Проведем из точки $x$ правую касательную к столу $\gamma$; пусть $v$ — точка касания. Тогда $y = T_\lambda x$ есть такая точка на продолжении луча $xv$, что $|vy| / |xv| = \lambda$. Нетрудно видеть, что $T_1$ есть имеющееся у нас преобразование внешнего биллиарда. В статье [15] исследуются внешние биллиарды со сжатием вокруг правильных n-угольников для n = 3, 4, 5, 6, 8, 12. Одним из основных результатов статьи является $\lambda$-стабильность всех периодических траекторий, т. е. для каждой периодической относительно $T_1$ орбиты существует $\varepsilon > 0$, т.ч. для всех $\lambda$, т.ч. $1-\varepsilon < \lambda < 1$ существует периодичная относительно $T_\lambda$ орбита с той же последовательностью вершин отражения; также доказана сходимость компонент при λ ↑ 1 в Хаусдорфовой топологии; в случаях же n = 3, 4, 6 для $T_\lambda$ при 0 < λ < 1 доказано существование конечного числа периодических орбит, к которым все остальные орбиты сходятся. Доказательство первого из этих утверждений для n = 5, 8, 12 базируется на существований полной ренормализационной схемы для «обычных» внешних биллиардов и рассмотрено лишь для случая n = 5; таким образом, чтобы завершить доказательства для n = 8 и 12, требуются полные ренормализационные схемы. Что касается случая правильного семиугольника, компьютерные эксперименты, проведенные Ёнгом, позволяют установить, что точки изображенной на рис.16 пятиугольной компоненты не являются λ-стабильными. Этот факт даёт возможность предполагать, что не имеет места быть и сходимость компонент в случае правильного семиугольника. Впрочем, отсюда не следует отсутствие ренормализационной схемы для других n.

Обсуждая внешние биллиарды вокруг многоугольников, нельзя также не упомянуть о «схеме развертки», аналоге метода развертки для внутренних биллиардов. Во внутреннем случае, при попадании шара в стенку, мы вместо отражения траектории от «стенки» отражаем сам стол относительно оси, содержащей «стенку»; во внешнем случае, мы оставляем на месте шар, а относительно точки касания центрально-симметрично отражаем сам стол (рис.18).

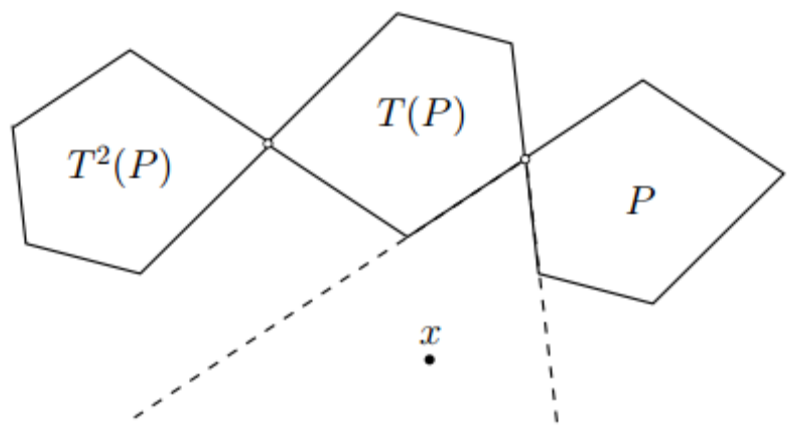

Рис. 18. Схема развертки

Заметим, что этот метод позволяет для каждой точки найти множество точек той же траектории. Действительно, отражение относительно фиксированной вершины может произойти в том и только в том случае, когда точка расположена внутри связанного с вершиной угла, изображенного на рис.18.. Следовательно, все точки, обладающие такой же траекторией, должны лежать в пересечении таких углов, причем такое условие является и достаточным. Такой подход даёт алгоритм, позволяющий для каждой периодической точки найти её периодическую компоненту абсолютно точно.

Закончим литературный обзор рассмотрением связи внешних биллиардов относительно многоугольников и символической динамикой. Перенумеруем вершины многоугольника числами от 1 до n; тогда орбита суть конечное или бесконечное слово, состоящее из букв алфавита {1, 2, 3, …, n}; множество всех конечных подслов всех слов, получающихся таким образом, составляют *язык* преобразования внешнего биллиарда вокруг заданного многоугольника; сложность $p = p(l)$ языка есть натуральнозначная функция, ставящая в соответствие числу l количество различных слов длины l в языке.

В [16] Гуткину и Табачникову удалось показать, что:

1) если стол внешнего биллиарда — правильный n-угольник, то существуют две положительные константы a, b такие, что для всех натуральных l верно

$$a*l \leq p(l) \leq b*l^{r+2},$$

где r есть функция Эйлера $\varphi(n)$;

2) если стол внешнего биллиарда — решеточный n-угольник, то существуют две положительные константы a, b такие, что для всех натуральных l верно

$$a*l^2 \leq p(l) \leq b*l^2.$$

В [17] для случаев n = 3, 4, 5, 6 удается найти точные описания языков внешних биллиардов относительно правильных n-угольников и найти их сложность; для случая n = 10 же доказывается биективность между языками для десяти- и пятиугольника; во всех случаях существует константа C > 0 (для каждого случая своя), т.ч. $p(l) \sim C*l^2$.

**Глава 2. Внешние биллиарды на многоугольниках: общие сведения**

Как уже было замечено во введении, при рассмотрении отображения внешнего биллиарда (будем называть его по-прежнему T) для многоугольников возникает проблема некорректности определения T на продолжениях сторон многоугольника; в этом случае будем говорить, что T не определено для таких точек. Таким образом, область определения T для n-угольника распадается на n областей $D_i$, в каждой из которых T есть центральная симметрия относительно соответствующей вершины $A_i$ многоугольника (рис. 19).

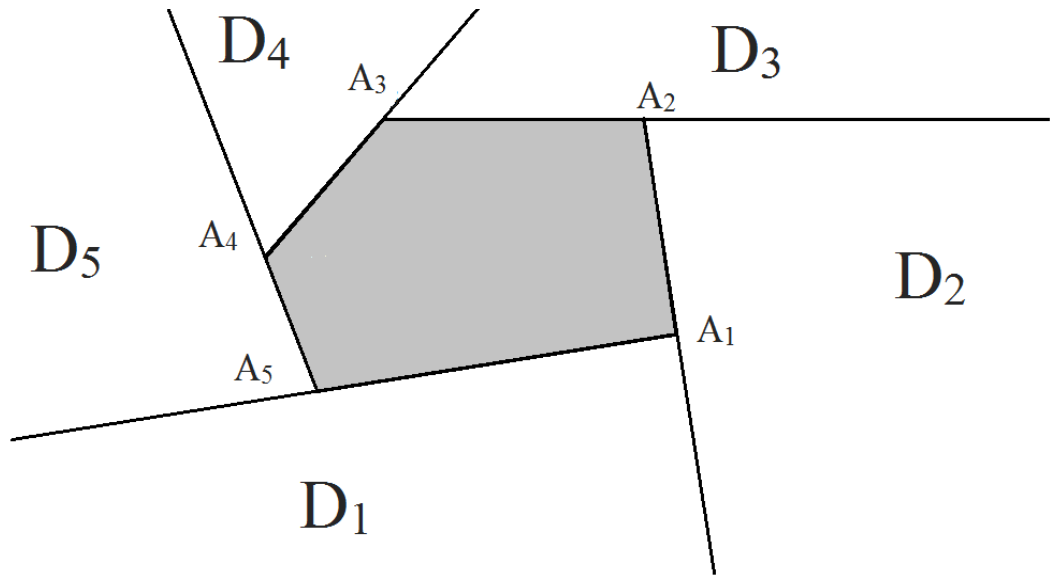

Рис.19 Области определения внешнего биллиарда для многоугольника

Напомним также, что в случае многоугольника точки вне стола можно разбить на следующие три типа: 1) точки с конечной траекторией (случай, когда на очередном шаге мы попадаем в точку, в которой T не определено); 2) точки с периодической траекторией; 3) точки с апериодической траекторией (апериодические, или бесконечные, точки).

Рассмотрим несколько простых свойств внешнего биллиарда:

1. Множества точек каждого типа инвариантны относительно T (с точностью до точек, в которых T не определено);

2. Пусть $B_i$ – множество точек, т.ч. $T^i$ определено, а $T^{i+1}$ нет. Тогда а) $B_{i+1} = T^{-1}(B_i \setminus \{$множество точек, в которых $T^{-1}$ не определено$\})$; б) множество точек первого типа есть $A_0 := \{$объединение всех $B_i\}$; в) $A_0$ есть счетное объединение открытых отрезков и лучей; очевидным следствием является тот факт, что размерность Хаусдорфа множества $A_0$ равна 1;

3. Если точка *x* имеет периодическую траекторию, то существует окрестность этой точки, состоящая целиком из периодических точек (например, можно взять (открытую) ε-окрестность точки *x*, где ε - минимальное из расстояний от точек траектории до границ соответствующих $D_v$ и половин попарных расстояний между точками траекторий); более того, если период *x* четен, то период точек проколотой окрестности будет совпадать с периодом *x*; в противном случае период точек окрестности будет вдвое больше;

3а. Если точка *x* имеет конечную траекторию, то существует открытый отрезок на плоскости, содержащий *x* и состоящий из точек с конечной траекторией (рассмотрим последнюю точку $x_{посл}$ траектории; эта точка лежит на продолжении стороны многоугольника. Очевидно, что можно взять достаточно малый отрезок этого луча, содержащий $x_{посл}$, и применить к этому отрезку обратное преобразование внешнего биллиарда несколько раз, вернув $x_{посл}$ в *x*, т.ч. отрезок при этом перешел в отрезок той же длины);

4. $A_0$ разбивает плоскость на компоненты, являющиеся открытыми выпуклыми фигурами (возможно, нулевой меры); каждая из этих компонент при преобразовании T переходит в равную компоненту;

5. Последовательность вершин стола, относительно которых происходит

отражение при построении траектории точки, зависит лишь от компоненты, в которой лежит эта точка; следовательно, понятия "типы 2 и 3" можно применять к компонентам;

6. Компонента нулевой меры не может быть компонентой 2-го типа (прямое следствие свойства 3);

7. Все точки периодической компоненты имеют один и тот же чётный период, кроме, возможно, одной точки; в этом случае компонента является центрально-симметричной, а «выколотой» точкой является центр симметрии компоненты, причём его период нечётен и равен половине периода остальных точек (при каждом применении T компонента переходит в центрально-симметричную; следовательно, в тот момент, когда она перейдёт сама в себя, либо все точки вернутся в своё первоначальное положение (первый случай), либо все точки перейдут в симметричные относительно центра симметрии компоненты (второй случай). Чётности периодов всех точек в обоих случаях очевидны; заметим, что центрально-симметричная компонента не обязана иметь выколотый в смысле периода центр – контрпример мы увидим далее);

8. Если траектория компоненты ненулевой меры ограничена, то она (траектория) периодична (ибо на ограниченном пространстве существует лишь ограниченное число равных фигур);

9. $T^2$ есть параллельный перенос вдоль некоторой стороны либо диагонали; длина вектора переноса есть удвоенная длина этой диагонали или стороны (как композиция центрально-симметричных отражений);

10. Компонента есть ограниченное множество. Действительно, предположим, что одна из компонент бесконечна. Применим аффинное преобразование так, чтобы вершина $A_1$ отражения T для этой компоненты (назовём её C) была самой «правой» вершиной многоугольника (рис.20).

Очевидно, что $T^2$ сдвинет компоненту «влево-вверх» либо строго «влево», причем длины сдвигов «влево» и «вверх» отделены от нуля. Заметим, что сдвиг строго «влево» при $T^2$ возможен лишь в случае, если вторая вершина отражения есть $A_n$, т.е. $T(C)$ лежит в $D_n$ (рис.21, 22).

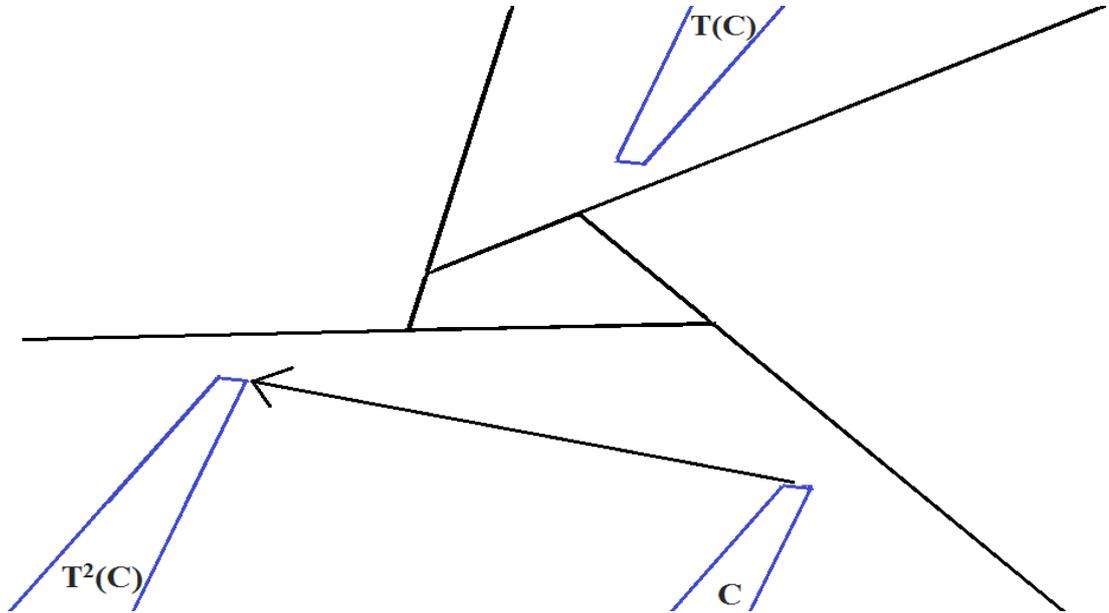

Рис. 20. Сдвиг бесконечной компоненты при $T^2$ «влево-вверх»

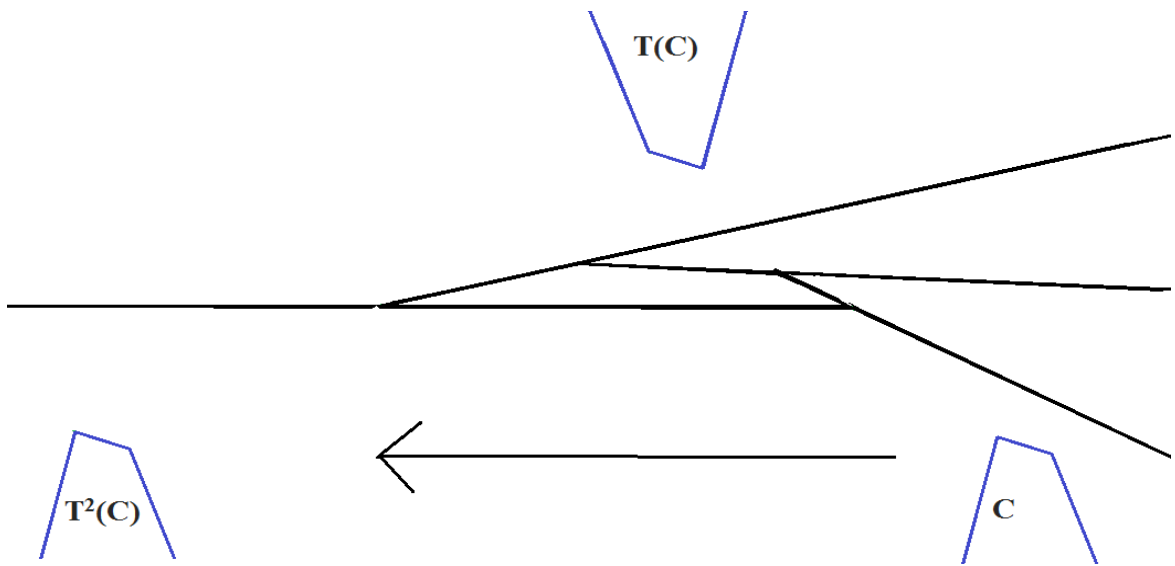

Рис. 21. Сдвиг бесконечной компоненты при $T^2$ «влево»

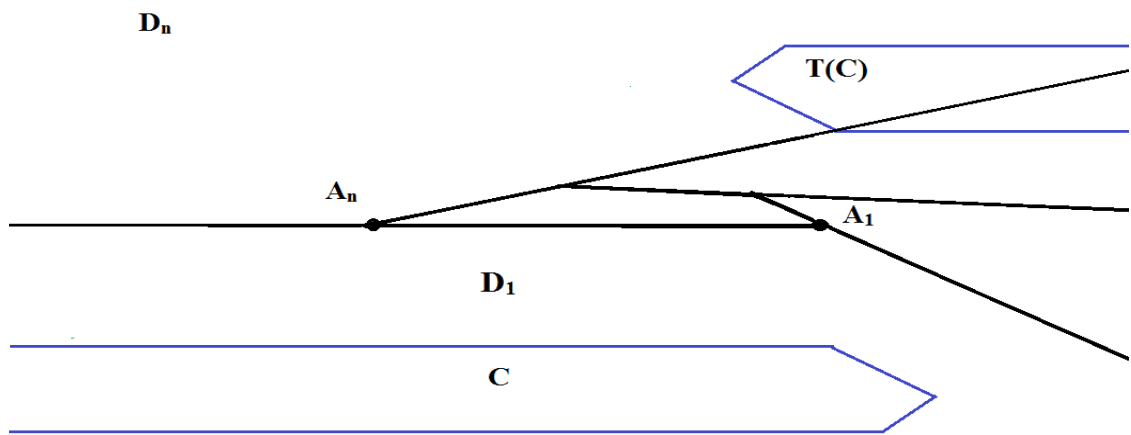

Рис. 22. «Бесконечная только влево» компонента

По мере применения $T^2(C)$ сдвигается влево, а $T(C)$ – вправо, причем на постоянную величину. Следовательно, через несколько итераций $T(C)$ покинет $D_n$, причем навсегда (ибо $T(C)$ сдвигается либо строго «вправо», либо «вправо-вниз»). Это означает, что $T^2$, начиная с некоторого момента, будет двигать C «влево-вверх». Отсюда, в свою очередь, следует, что рано или поздно C покинет $D_1$ и переместится в $D_n$. Т.к. $D_n$ ограничена снизу, то и C не может быть неограниченной снизу. Т.к. C в исходном состоянии ограничена лучом $A_1 A_n$ сверху и $A_2 A_1$ справа, то она (C) неограничена слева. Применим к C последовательность операций T, переводящие C в $D_n$, после чего «повернем» картинку так, чтобы ребро $D_{n-1} D_n$ стало «нижним»; тогда и C окажется «неограниченной снизу»; повторение рассуждений доказывает свойство.

Однако для этого свойства существует ещё одно интересное доказательство, использующее «схему развертки». Действительно, будем вместо отражения относительно точки касания точки-шара отражать стол; тогда вся компонента содержится в углу между лучами — продолжениями последовательных сторон (см. рис. 18). Более того, компонента является пересечением всех таких углов. Соединим последовательные точки отражения отрезками; очевидно, что каждый вышеописанный угол лежит в полуплоскости, задаваемой прямой — продолжением соединяющего вершину угла с предыдущей точкой отражения отрезка; так как точка обладает бесконечной (возможно, периодичной) траекторией, то стол, а вместе с ним и последовательность «обойдет» точку хотя бы раз; следовательно, пересечение полуплоскостей, задаваемых ребрами, а также и целиком лежащая в этом пересечении компонента, суть ограниченные множества, QED.

11. Каждая компонента представляет собой либо (открытый) выпуклый многоугольник, стороны которого параллельны сторонам многоугольника, либо отрезок, параллельный одной из сторон многоугольника, либо точку (очевидное следствие применения схемы развертки).

Относительно этого свойства отметим, что автор не встречался со случаем компоненты отрезка. Теоретически, ни одно из свойств напрямую не противоречит ситуации, когда компонента ненулевой площади имеет уходящую в бесконечность траекторию; однако в силу Теоремы 1 выполнено следующее свойство:

12. Траектории внешнего биллиарда вне решеточных и правильных многоугольников ограничены.

13. Все компоненты для решеточного многоугольника суть невырожденные многоугольники (ибо верно УТВ.: расстояния между параллельными прямыми, содержащими лучи либо отрезки множества $A_0$, отделены от нуля. Докажем УТВ. для лучей, параллельных ребру $A_nA_1$ произвольного выпуклого решеточного многоугольника. Переведем $A_nA_1$ в горизонтальный отрезок так, чтобы вершины многоугольника остались целочисленными. При таком преобразовании расстояния между всеми парами прямых, содержащих параллельные $A_nA_1$ лучи/отрезки, а) домножились на одну и ту же константу; б) превратились в натуральные числа, т.е. стали больше либо равны 1 ( в силу коммутируемости T и аффинного преобразования). Следовательно, УТВ.

доказано, а с ним и свойство).

14. Все компоненты решеточного многоугольника есть невырожденные компоненты типа 2 (прямое следствие свойств 8, 12, 13).

Вооруженные таким мощным багажом знаний о внешних биллиардах, мы можем перейти к изучению конкретных примеров многоугольников.

**Глава 3. Внешний биллиард вне квадрата**

По-видимому, квадрат является если не самым простым для исследования внешнебиллиардным столом, то по крайней мере самым простым многоугольным столом. Попытки нарисовать $A_0$ вручную, равно как и компьютерные эксперименты, дают возможность предполагать, что множеством $A_0$ в случае квадрата с вершинами (0, 0), (0, 1), (1, 1), (1, 0) является целочисленная сетка (т.е. набор прямых вида $x = C$ и $y = D$, где $C$ и $D$ – целочисленные константы) (рис.23).

Прямым следствием этой гипотезы является тот факт, что получившиеся «квадратики» должны стать компонентами. Попробуем доказать такую гипотезу.

Факт о том, что при преобразовании $T$ квадратик переходит в квадратик, вполне очевиден. Однако напрямую это означает лишь то, что каждый квадратик принадлежит одной компоненте, но не то, что каждый квадратик есть отдельная компонента. Но из рис.24 очевиден следующий факт: при применении $T$ не изменяется «манхэттенское расстояние» от квадратика до стола.

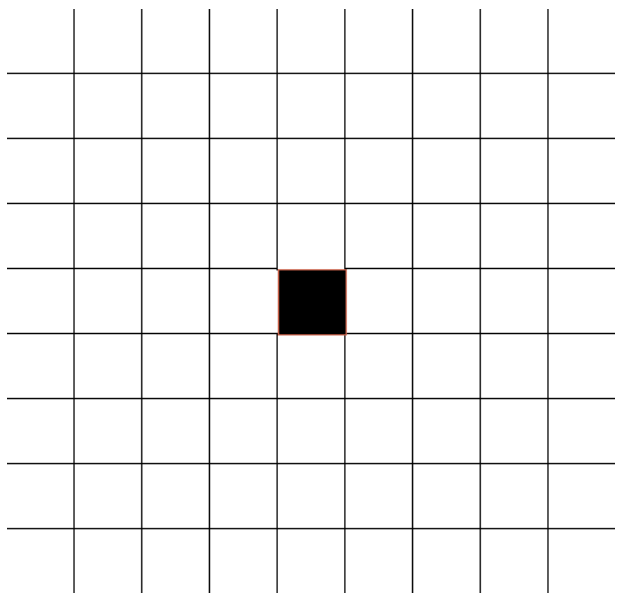

Рис. 23: Потенциальное множество точек первого типа

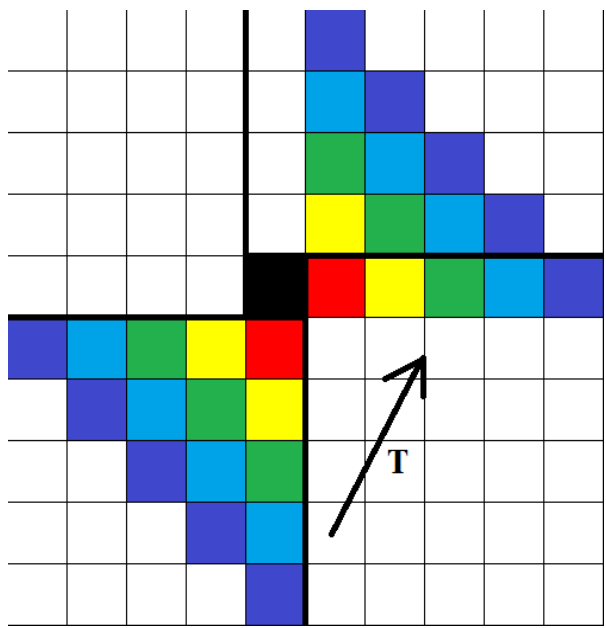

Рис. 24. Иллюстрация инвариантности «манхэттенского» расстояния до стола

Итак, каждый квадратик движется строго по «ожерелью» из 4d квадратиков, находящихся на (манхэттенском) расстоянии d. Если раскрасить все эти квадратики в шахматном порядке (рис. 25), то становится видно, что 1) преобразование T меняет цвет квадратика; 2) преобразование $T^2$ перемещает квадратик по ожерелью на два квадратика по часовой стрелке.

Последние два факта позволяют заключить, что период каждого квадратика есть ровно 4d, а так как любые два соседних квадратика имеют разные периоды, то граница между ними состоит из точек лишь первого типа (ибо в какой-то момент квадратики будут отражаться от разных точек), что и приводит нас к итоговой картинке, ранее предсказанной компьютерными экспериментами, изображенными на рис.26.

Аналогичным образом можно провести анализ для случаев правильных треугольника и шестиугольника. Здесь можно увидеть и компоненты с выколотым по периоду центром, и нецентрально-симметричные компоненты. В случае шестиугольника можно видеть, что шестиугольные компоненты одного периода делятся на две орбиты, по 3*level шестиугольников в каждом, а треугольные есть единая орбита из 12*level-6 треугольников. В треугольном же случае мы имеем орбиты из 6*level-3 шестиугольников и 12*level треугольников. Компьютерные эксперименты в этих случаях изображены на рис.26.

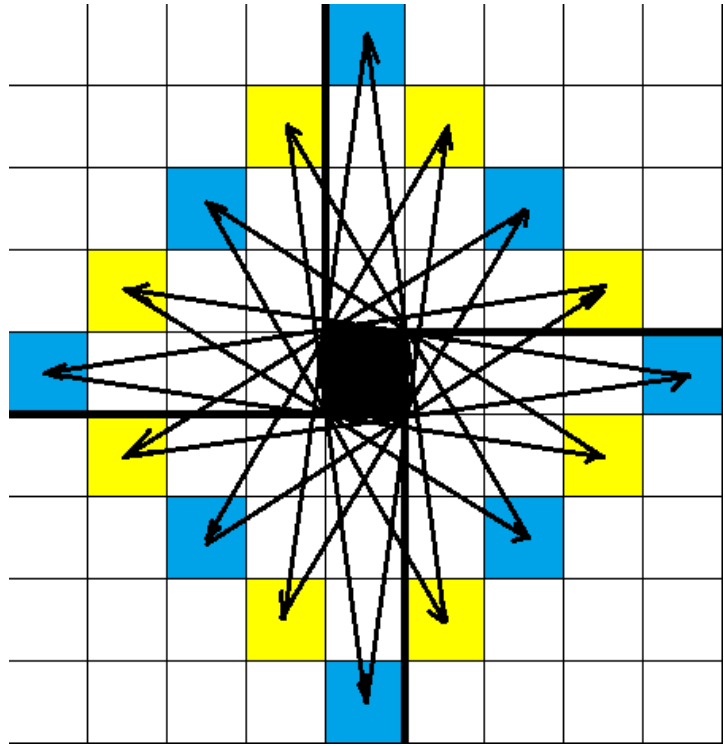

Рис. 24. Траектория одного из квадратиков

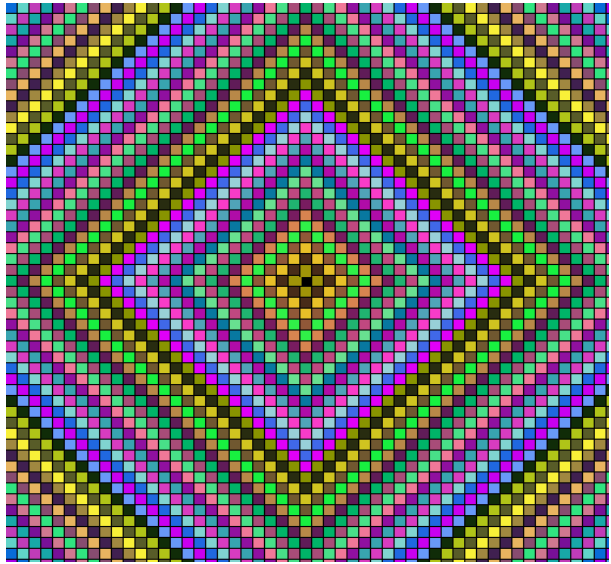

Рис. 25: Зоопарк точек внешнего биллиарда вне четырехугольника

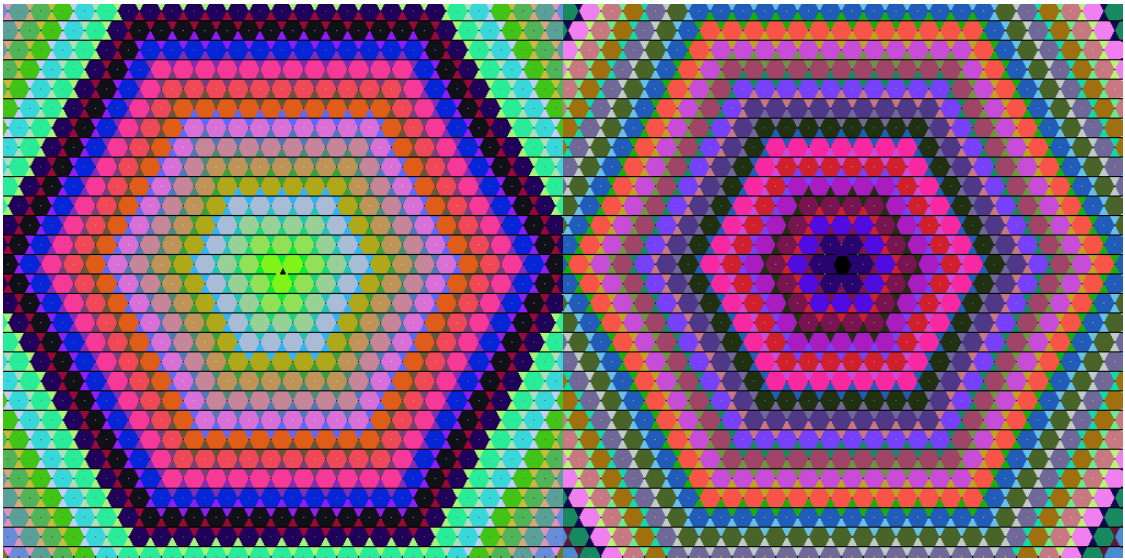

Рис. 26. Внешний биллиард вне правильных шестиугольника и треугольника

## Глава 4. Внешний биллиард вне правильного восьмиугольника

Простейшим случаем правильного нерешеточного многоугольника является правильный пятиугольник. Этот случай был подробно исследован в, например, [1]. Мы же проведем аналогичное исследование для правильного восьмиугольника.

Рассмотрим следующую изображенную на рис.27 инвариантную относительно T компоненту :

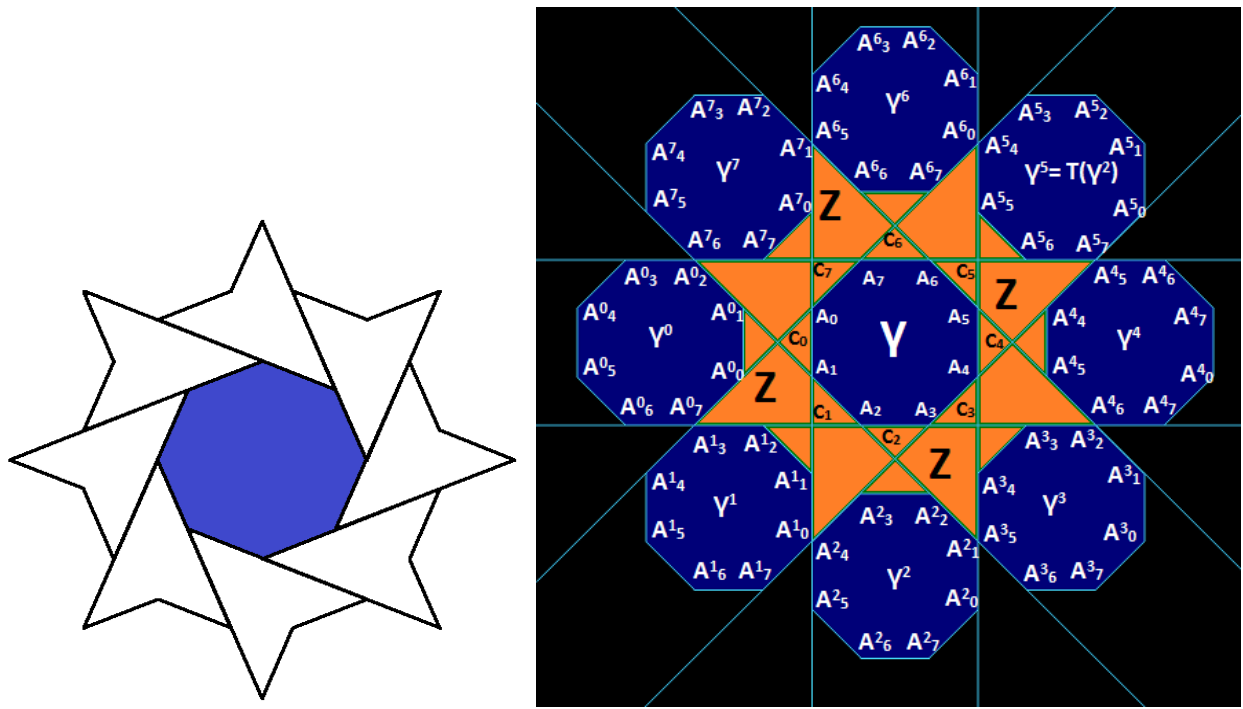

Рис.27. Инвариантная относительно T компонента Z

Разделим ее на 8 равных частей, как показано на рис.27, и отождествим их относительно поворота на 45n градусов; будем понимать под T производное отображение на получившейся фигуре. Как видно из рис. 28., преобразование T поворачивает треугольник POQ на 135 градусов относительно точки U, четырехугольник KPQR – на 90 градусов вокруг V, а треугольник LRM – на 45 градусов вокруг точки W (во всех случаях поворачиваем против часовой стрелки). Назовем эти повороты u, v и w соответственно.

Заметим, что четырехугольник OKLM можно разбить на «вписанный» в него правильный восьмиугольник с центром в точке V и три равные фигуры, подобные OKLM. Каждую из этих фигур можно разбить рекурсивно (рис.29).

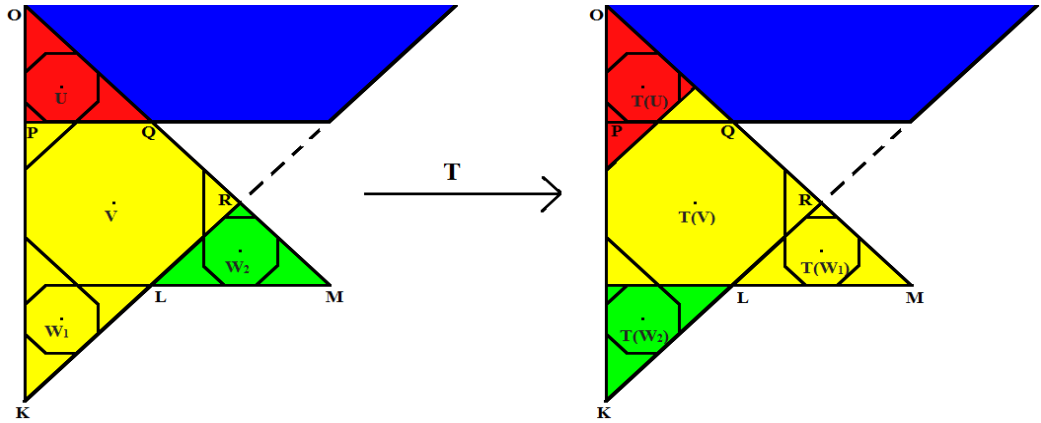

Рис.28 Модифицированное преобразование внешнего биллиарда

Заметим, что центральный восьмиугольник инвариантен относительно T. Верно также и то, что множество точек, лежащих в восьмиугольниках одного размера, также инвариантно относительно T. Не будем проводить доказательство этого (очевидного) факта; вместо этого сосредоточимся на дальнейшем анализе. Введем преобразование Γ, являющееся сжатием с центром в т. O и переводящим OKLM в OK'L'M' (рис.30).

Из рис.21 очевидно получается следующая

**Лемма 1:** $\Gamma T x = T^k \Gamma x$, где k = 15 для x-ов треугольника OPQ, 9 для x-ов четырехугольника KPQR и 3 для x-ов треугольника LRM. Более точно,

$\Gamma u(x) = uvvwvwvwvwvvu\Gamma x$, $\Gamma v(x) = uvvwvwvvu\Gamma(x)$, $\Gamma w(x) = uuu\Gamma(x)$

Для дальнейшего анализа введем понятие ранга. Рангом точки x назовем максимальное n такое, что $\Gamma^{-n} x$ еще лежит в четырехугольнике OKLM, а рангом орбиты – максимум среди рангов всех ее точек.

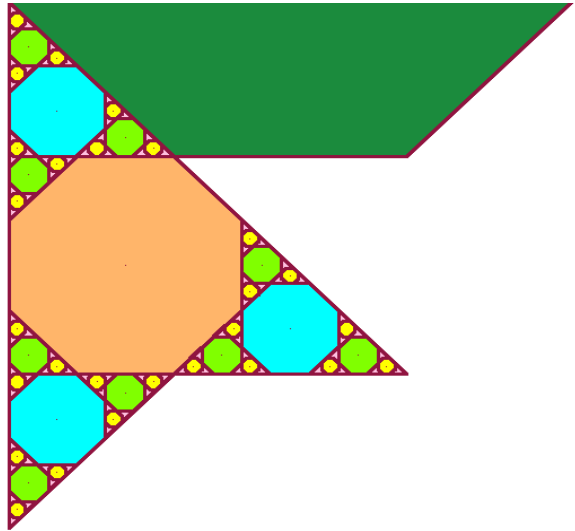

Рис.29. Самоподобие фигуры

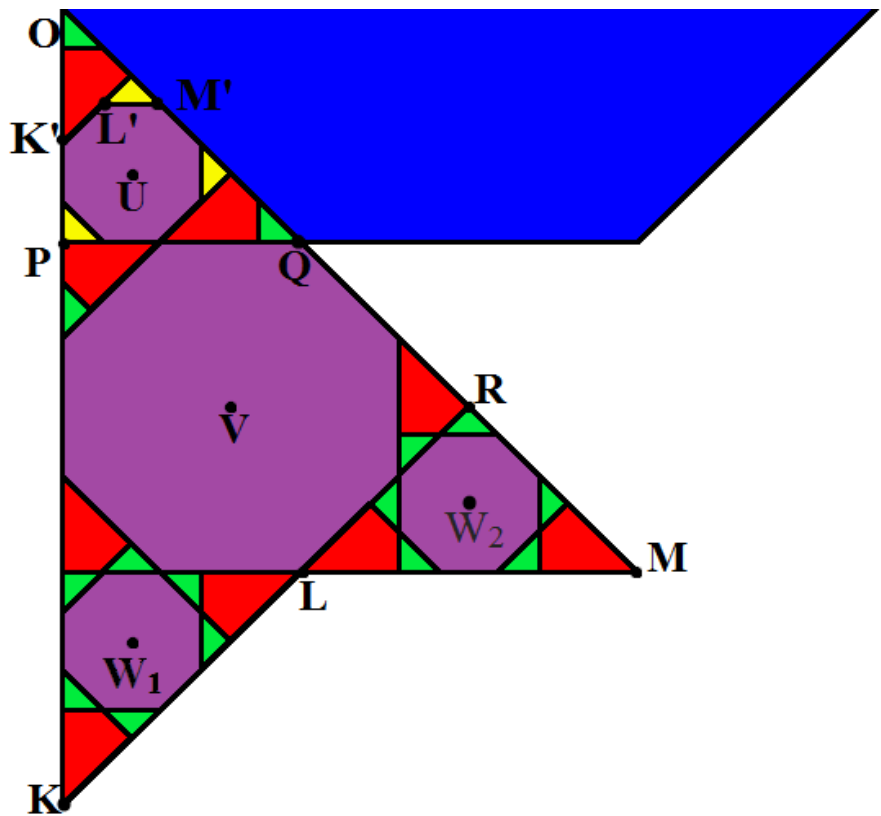

Рис.30. Доказательство леммы 1

**Лемма 2:** любая траектория ранга n > 0 может быть получена из траектории ранга n-1 путем подстановки по правилу Г из леммы 1.

**Доказательство:** рассмотрим орбиту ранга n, и пусть x – точка ранга n этой орбиты, а y = Г$^{-1}$x. Тогда по лемме 1, ГT$^k$y = T$^{f(k)}$x, где f(k) – некая возрастающая функция. Т.к. ни одна из точек ГT$^k$y не имеет ранга > n, то ранг траектории y не превышает n-1, а т.к. ранг y есть n-1, то лемма доказана.

Лемма 2 означает, что любая периодическая траектория ранга n получается из траектории ранга 0 путем n применений операции Г. Заметим также, что при применении Г две соседние точки x1 и x2 орбиты превращаются в две точки орбиты ранга выше, между которыми появляется несколько (2, 8 или 14, если быть точнее) вершин ранга 0 (очевидно из рис.30). Это дает нам возможность посчитать размер любой орбиты. Сделаем это следующим образом. Пусть в текущей периодической орбите для получения следующего элемента $a_k$ раз применяется оператор u, $b_k$ раз применяется оператор v, $c_k$ раз применяется оператор w (k – ранг траектории). Тогда из леммы 1 несложно увидеть, что после применения к орбите Г получаем:

$$a_{k+1} = 2a_k + 2b_k + 3c_k,\ b_{k+1} = 8a_k + 5b_k,\ c_{k+1} = 5a_k + 2b_k.$$

Разрешая эту систему, получаем:

$$a_k = (1 + 4(-3)^k + 3*9^k)a_0 + (-2 + 2*9^k)b_0 + (3 - 4(-3)^k + 9^k)c_0,$$

$$b_k = (-2 - 4(-3)^k + 6*9^k)a_0 + (4 + 4*9^k)b_0 + (-6 + 4(-3)^k + 2*9^k)c_0,$$

$$c_k = (1 - 4(-3)^k + 3*9^k)a_0 + (-2 + 2*9^k)b_0 + (3 + 4(-3)^k + 9^k)c_0,$$

а величина орбиты ранга k есть

$$a_k + b_k + c_k = (1.5*9^n - 0.5*(-3)^n)a_0 + 9^{n*}b_0 + (1.5*9^n + 0.5*(-3)^n)c_0.$$

Остается лишь рассмотреть, какие траектории ранга 0 имеются в наличии; все они перечислены в таблице 1 (окрестностью в каждом случае является соответствующий восьмиугольник; траектории ранга 0 на рис.30 раскрашены в фиолетовый цвет).

Таким образом, множество точек, имеющих периодические орбиты n-го уровня, есть, как несложно видеть, набор из $4*9^n$ восьмиугольников, причем размер орбиты экспоненциально растет с ростом ранга.

Нетрудно видеть, что существуют точки четырехугольника OKLM, не лежащие внутри одного из таких восьмиугольников или на их границах; среди них и нужно искать точки с бесконечной апериодической траекторией.

Таблица 1. Траектории ранга 0 и их характеристики

| Положение стартовой точки | Маршрут | $a_0$ | $b_0$ | $c_0$ | Период соответствующей траектории ранга n |
|---|---|---|---|---|---|
| Точка V | v | 0 | 1 | 0 | $9^n$ |
| Окрестность V | vvvv | 0 | 4 | 0 | $4*9^n$ |
| Точка U | U | 1 | 0 | 0 | $1.5*9^n - 0.5*(-3)^n$ |
| Окрестность U | Uuuuuuuu | 8 | 0 | 0 | $12*9^n - 4*(-3)^n$ |
| Точка $W_1$ | Vw | 0 | 1 | 1 | $1.5*9^n + 0.5*(-3)^n$ |
| Окрестность $W_1$ | $(vw)^8$ | 0 | 8 | 8 | $12*9^n + 4*(-3)^n$ |

Рассмотрим одну из таких точек, являющуюся пределом самоподобной последовательности точек, начинающейся с точек с0, с1, с2 (см. рис. 31), и продолжающейся аналогичным образом в подобном OKLM четырехугольнике, находящемся внутри треугольника с0с1с2 (т.е. есть перевести с помощью сжатия OKLM в такой четырехугольник, то с0, с1 и с2 перейдут в с3, с4 и с5, соответственно; дальнейшим сжатием получим с6, с7 и с8 и т.д.); назовем эту точку C. Эта точка не может обладать периодической орбитой (ибо иначе C должна обладать проколотой окрестностью точек одинакового периода, что в данном случае неверно, ибо периоды компонент cN растут бесконечно), равно как и конечной орбитой (ибо в таком случае должен существовать отрезок, содержащий строго внутри себя C, содержащий точки с орбитами того же размера и типа, причем содержащийся внутри всех треугольников с0с1с2, с3с4с5, …, что также не выполняется; лежать же в треугольнике сей отрезок должен потому, что в противном случае одна из точек пересечения отрезка с границей треугольника будет периодической точкой); следовательно, эта точка является точкой с бесконечной апериодической траекторией.

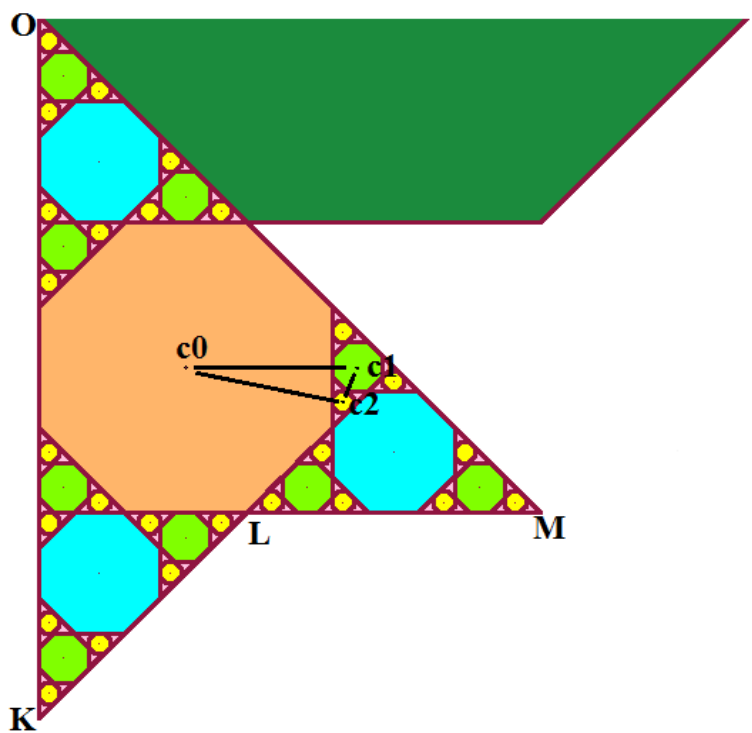

Рис.31. Периодические точки T.

Таким образом, мы показали, что существует точка с бесконечной апериодической траекторией для преобразования внешнего биллиарда относительно правильного восьмиугольника. С другой стороны, несложно показать, что периодические компоненты образуют в четырехугольнике OKLM множество полной меры. Но что же происходит за пределами рассматриваемой инвариантной компоненты?

Чтобы понять это, рассмотрим снова рисунок 27. До текущего момента, мы исследовали индуцированное преобразование T лишь в многоугольнике $A_1A^2_4A^2_3A^2_2$, который мы переобозначили как OKLM; теперь будем рассматривать T для бесконечного угла $A^2_4A_1A^2_1$. Заметим, что такой угол состоит из фигуры $A^2_5A^2_4A_1A^2_1A^2_0$ и угла $A^2_0A^2_1A^3_6$. Пусть T'' – преобразование первого возвращения T для угла $A^2_0A^2_1A^3_6$, а S – параллельный перенос, переводящий точку $A^1$ в точку $A^2_1$, а угол $A^2_4A_1A^2_1$ – в угол $A^2_0A^2_1A^3_6$.

**Лемма 3:** Пусть точка p лежит строго внутри угла $A^2_4A_1A^2_1$. Тогда:

1. T(p) определено, если и только если T''(S(p)) определено;

2. Если T(p) определено, то T''(S(p)) = S(T(p)).

Доказательство леммы 3 легко получить, просто рассмотрев траектории первого возвращения $A^2_0A^2_1A^3_6$ и то, на какие фигуры сей угол разбивается при последовательном применении индуцированного преобразования T. Более того, можно заметить, что траектории возвращения покрывают весь угол $A^2_4A_1A^2_1$, за исключением активно исследованного нами четырехугольника $A_1A^2_4A^2_3A^2_2$ и восьмиугольника $\gamma^2$. Отсюда, с помощью технологии, использованной для описания периодических компонент, можно доказать следующие леммы.

**Лемма 4:** для любой периодической(апериодической, граничной) точки q угла $A^2_4A_1A^2_1$ существуют целые неотрицательные числа n, k и точка p, лежащая в четырехугольнике $A_1A^2_4A^2_3A^2_2$ или в восьмиугольнике $\gamma^2$, т.ч. $q = T^k(S^n(p))$.

**Лемма 4а**(прямое следствие леммы 4) множество периодических точек как внутри угла $A^2_4A_1A^2_1$ относительно индуцированного T, так и вне стола γ для «неиндуцированного» T, есть множество полной меры.

Отметим, что вопрос о полноте периодических точек для внешнего биллиарда вне правильного n-угольника остается в общем случае открытым; известны лишь положительные решения вопроса для случаев n = 3, 4, 6 (отсутствие апериодических точек) и n = 5 (Табачников).

**Лемма 4б**(все также прямое следствие леммы 4): любая периодическая фигура для внешнего биллиарда вне правильного многоугольника γ равна некоторой периодической фигуре внутри четырехугольника OKLM = $A_1A^2_4A^2_3A^2_2$. В частности, все периодические компоненты суть правильные восьмиугольники.

Отсюда же, проанализировав подстановки, индуцируемые Г и S, можно вычислить и периоды фигур как относительно индуцированного преобразования T, так и для «неиндуцированного» T, подсчитав количество поворотов типа u, v, w и т.д.. Получаем следующий результат.

**Лемма 5.** Множество различных периодов для внешнего биллиарда вне правильного восьмиугольника есть объединение следующих четырех множеств:

- $\{12*9^n-4*(-3)^n, 12*9^n+4*(-3)^n, 4*9^n \mid n \in \mathbb{Z}, n \geq 0\}$;

- $\{8, 8k, 8k*9^n \mid n, k \in \mathbb{Z}, n \geq 0, k \geq 2\}$;

- $\{24k*9^n-4*(-3)^n+12*9^n \mid n, k \in \mathbb{Z}, n \geq 0, k \geq 2\}$;

- $\{24k*9^n+4*(-3)^n-4*9^n | \ n, k \in \mathbb{Z}, n \geqslant 0, k \geqslant 2\}$.

## Глава 5. Внешний биллиард вне правильного двенадцатиугольника

Исследуем на предмет нахождения бесконечной апериодической траектории внешний биллиард относительно двенадцатиугольника. Для начала выделим первую T-инвариантную компоненту, отождествим точки в этой компоненте, получающиеся друг из друга поворотом на 30 градусов, и рассмотрим фигуру, на которой определяется производное преобразование T (рис.32).

Мы видим, как «самолётик» для восьмиугольника превратился в «ракету» (назовём её A) для двенадцатиугольника; преобразование T делит «ракету» A на пять частей, или «зон T», как мы их иногда будем называть; верхняя часть является равнобедренным треугольником, три средние похожи на «воздушных змеев» Шварца, и действительно являются ими; пятая же часть оказывается невыпуклой, что является первым важным отличием случая двенадцатиугольника от случаев пяти- и восьмиугольника, где область определения T делилась лишь на выпуклые фигуры. Перенумеруем зоны T от 0 до 4 в порядке удаления от вершины «ракеты»; тогда для зоны T номер $i$ преобразование T является поворотом на $150 - 30*i$ градусов.

Заметим, что образы первых четырех зон T при применении T имеют пересечения с прообразами. Последовательно применяя преобразование T и пересекая образ с изначальной зоной, для каждой из этих четырех зон можно получить «инвариантную фигуру зоны $i$ преобразования T» - часть зоны, инвариантную относительно T. Такие инвариантные фигуры изображены на рис. 33.

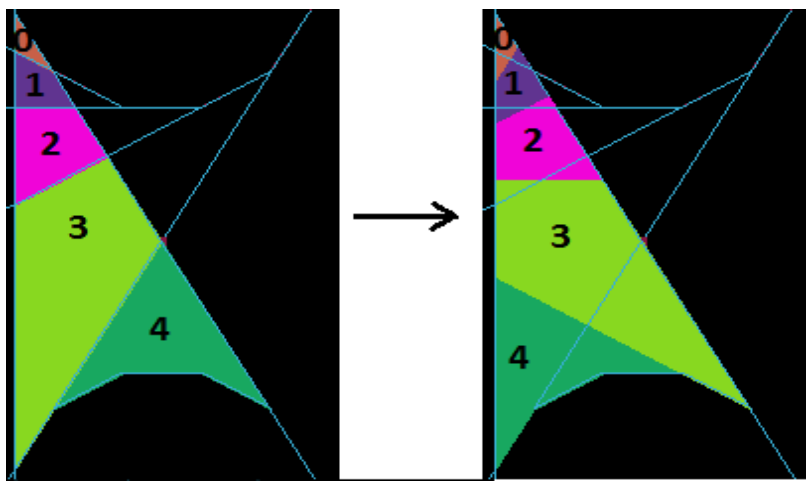

Рис.32. Ограничение преобразования внешнего биллиарда для правильного двенадцатиугольника

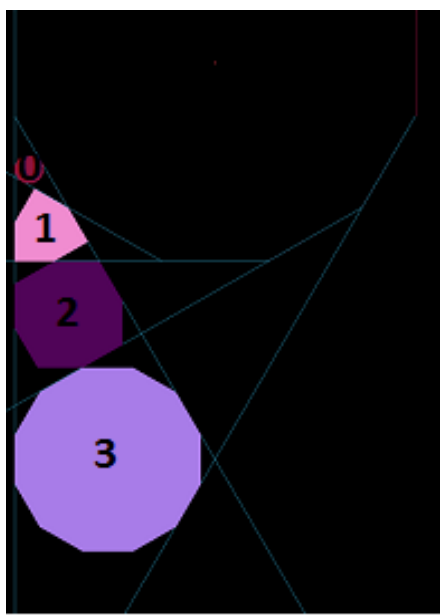

Рис.33: Инвариантные фигуры зон 0, 1, 2, 3

Здесь обнаруживается еще одно принципиальное отличие случая двенадцатиугольника. А именно: в случае пятиугольника все периодические компоненты являются правильными пяти- и десятиугольниками; в случае восьмиугольника все компоненты суть правильные восьмиугольники; в случае же внешних биллиардов относительно правильного двенадцатиугольника инвариантные многоугольники могут быть неправильными! Более точно, только две из инвариантных фигур являются правильными двенадцатиугольниками. Инвариантная фигура номер один (нумерация аналогична оной для зон T') равносторонним, но неправильным шестиугольником с углами по 90 и 150 градусов, инвариантным относительно поворота на 120 градусов; инвариантная же фигура номер два оказывается неправильным равносторонним восьмиугольником с углами по 150 и 120 градусов, инвариантным относительно поворота на 90 градусов.

По аналогии со случаем восьмиугольника, оказывается целесообразным ввести преобразование $Г_\lambda$, сжимающее большую ракету A в подобную ей маленькую $A_{мал}$, располагающуюся над инвариантной фигурой зоны 0 и изображенную на рис. 34; коэффициент сжатия $\lambda$ равен $7 - 4*sqrt(3) \approx 0.072$.

Поймём, каким образом могут выглядеть периодические компоненты компоненты ранга 0 (относительно $Г_\lambda$). По аналогии с восьмиугольником, такими компонентами являются уже известные нам инвариантные фигуры, а также фигуры, получающиеся из изначальных путем поворота на 150 градусов вокруг центра инвариантной фигуры зоны номер 3; все эти фигуры и их нумерация изображены на рис.35.

Отметим, что уже на этом этапе мы видим периодические траектории компонент с периодами 3 и 4. Можно заметить, что существуют еще две достаточно большие компоненты ранга 0; эти компоненты (и их траектории) изображены на рис.36.

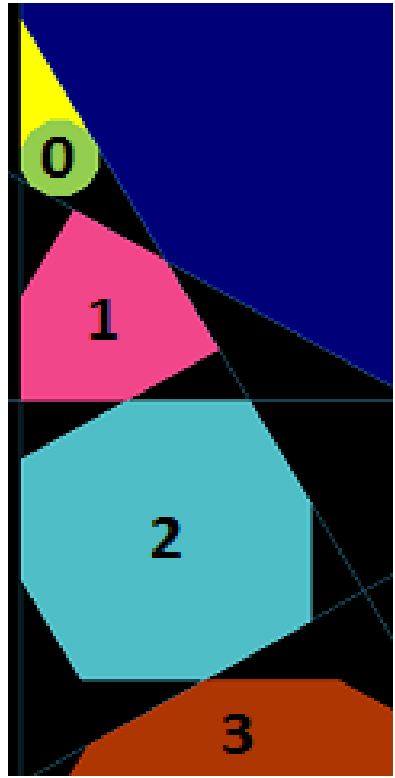

Рис. 34. Маленькая ракета $A_{мал}$ (выделена желтым цветом)

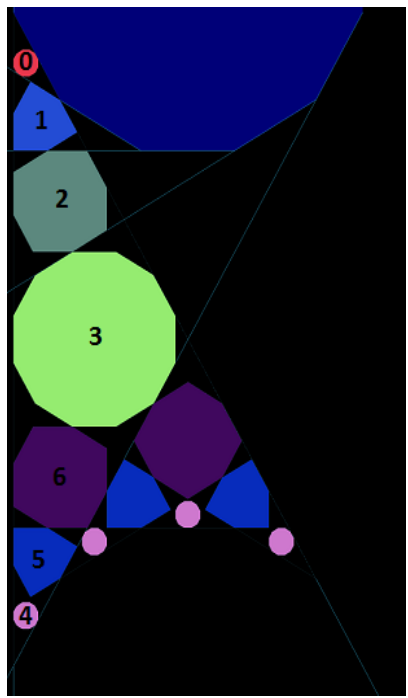

Рис. 35. Первые 7 периодических компонент

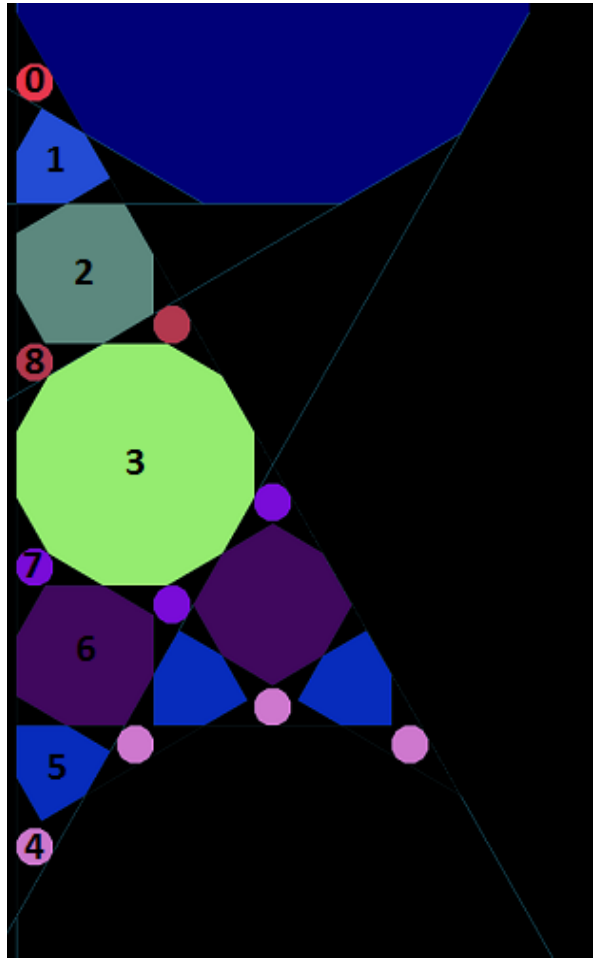

Рис. 36. Первые девять периодических компонент.

Возникает гипотеза о том, что только эти компоненты являются компонентами ранга 0. Проверим её. Для этого сожмем имеющиеся компоненты в λ (точнее, в 1/λ) раз и рассмотрим траектории получившихся фигур на рис. 37. Из рисунка видно, что какие-то компоненты ранга 0 мы упустили; например, таковой компонентой может являться двенадцатиугольник, равный по размеру сжатой инвариантной фигуре зоны 0, на "кончике хвоста" большой ракеты A; самым же интригующим здесь является наличие "солнечной короны" вокруг двенадцатиугольников – членов траектории компоненты номер 4; по-видимому, внутри набора этих корон существуют и другие периодические траектории ранга 0; не будем углубляться далее в этот вопрос.

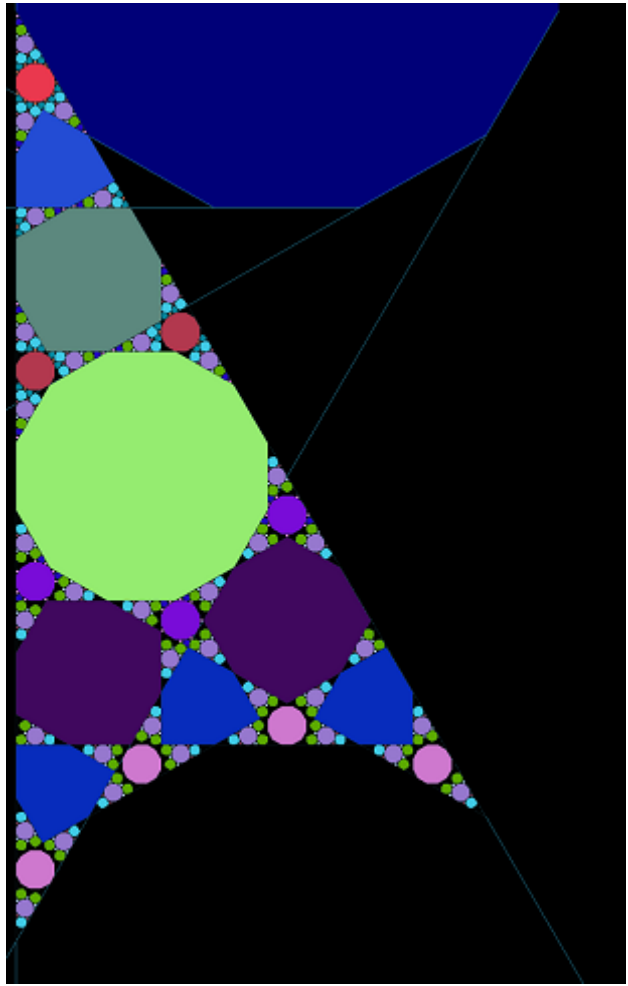

Рис. 37. Периодические траектории ранга 1

Хотя нам не удалось получить все траектории первого уровня, мы видим определенное самоподобие среди таких траекторий. В случае восьмиугольника найденная бесконечная апериодическая точка являлась пределом некоторой спиралеобразной последовательности компонент; в случае же двенадцатиугольника последовательность может выглядеть следующим образом. На рис.38 изображены три периодические компоненты, окружающие «меньшую среднюю ракету», подобную исходной; назовем эту ракету X, а компоненты, в порядке убывания размера – $C_0$, $C_1$, $C_2$. Пусть $\Gamma_X$ – преобразование, переводящее "среднюю ракету" $A_{ср}$ в "меньшую среднюю" X; данные из предыдущих рисунков дают возможность предположить, что три выделенные компоненты при применении $\Gamma_X$ перейдут в три многоугольника, также являющихся периодическими компонентами; более того, при последовательном применении Г$_X$ таким образом, по-видимому, можно получить все новые периодические компоненты с возрастающим периодом.

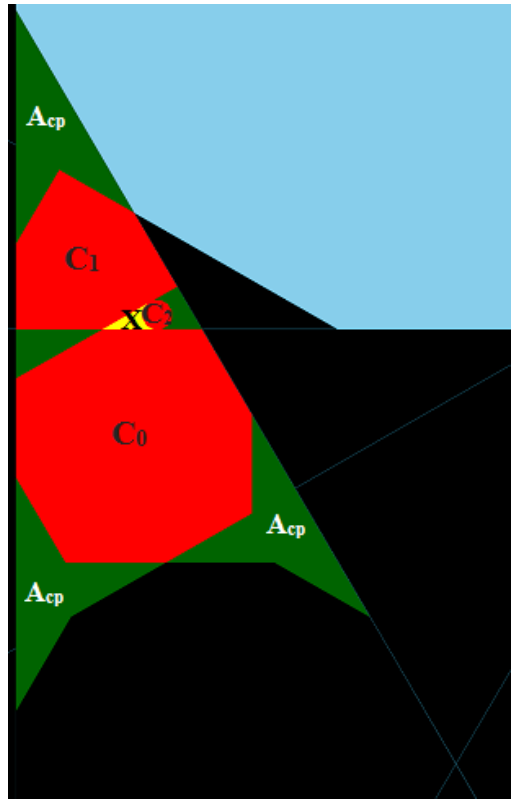

Рис.38. «Меньшая средняя» ракета X

Приступим к поиску бесконечной апериодической траектории. Для этого сформулируем алгоритм, который уже помог нам обнаружить апериодическую траекторию для правильного восьмиугольника. В уже имеющихся терминах, для двенадцатиугольника алгоритм выглядит следующим образом:

1) Найти в «ракете» A подобную ей меньшую «ракету» B и преобразование подобия Г, что для всех точек $x$ внутри «ракеты» A, где определено преобразование T, выполнено $ГTx = T^{k(x)}Гx$; заметим также, что для точки $y = Гx$ преобразование $T^{k(x)}y$ является преобразованием первого возвращения (first return map) относительно преобразования T. Таким образом, преобразование T' первого возвращения для маленькой «ракеты» должно быть устроено так же, как и преобразование Т, т. е. должно выполняться свойство ГТ = Т'Г;

2) Ввести относительно преобразования Г ранг орбиты; найти все периодические орбиты ранга 0;

3) Показать, что Г инвариантно относительно каждого из трех типов точек; вывести отсюда, что все периодические орбиты можно получить из орбит ранга 0 с помощью «подстановочного» преобразования Г; понять, как устроены периодические компоненты и их периоды;

4) Найти удачную последовательность сходящихся к некоторой компонент возрастающего периода; показать, что такая точка не может быть ни периодической, не обладать конечной траекторией; вывести отсюда, что такая предельная точка обладает бесконечной апериодической траекторией.

Попробуем реализовать такой алгоритм. По аналогии со случаями пяти- и восьмиугольников «хорошим» должно оказаться преобразование $Г_\lambda$, рассмотренное нами выше. Однако, компьютерные эксперименты показывают, что преобразование первого возвращения для «маленькой ракеты» выглядит совершенно по-другому.

Как видно из рис.39, преобразование первого возвращения делит «маленькую ракету» $А_{мал}$ на десять многоугольников; количества углов итераций до первого возвращения этих многоугольников описаны в таблице 2. Разумеется, хорошего преобразования Г в данном случае не существует.

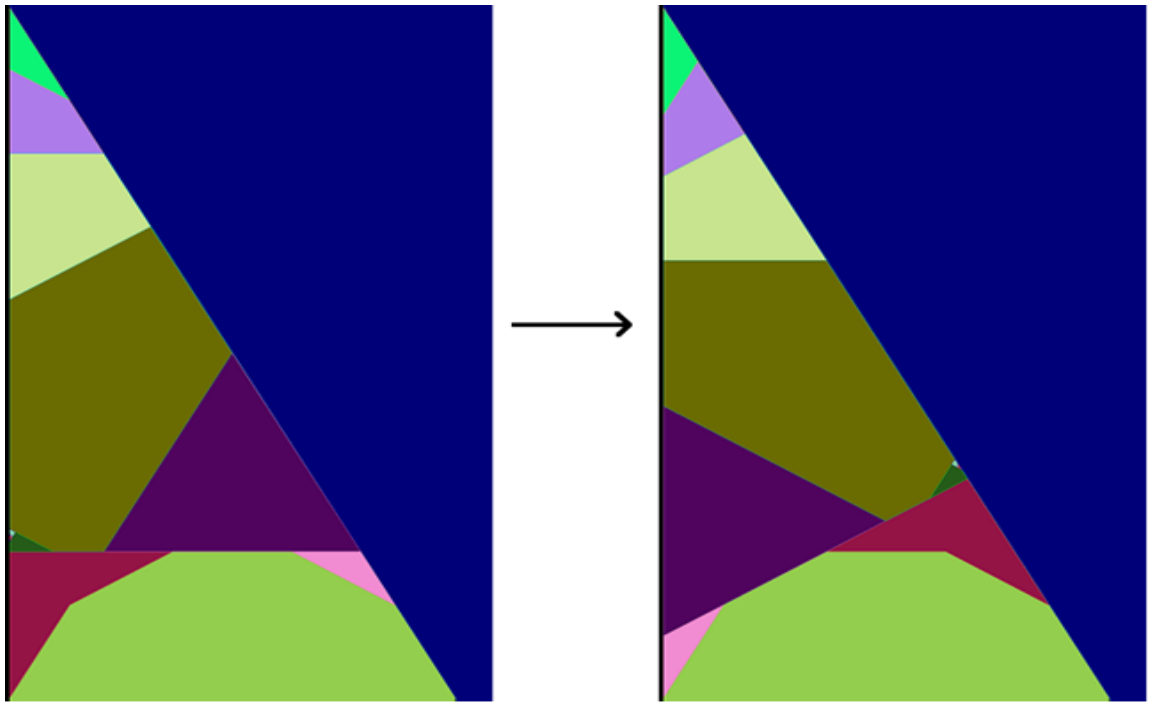

Рис. 39. Преобразование первого возвращения для $А_{мал}$

**Таблица 2**. Преобразование первого возвращения для $А_{мал}$

|  | 1 | 2 | 3 | 4 | 5 | 6 | 7 | 8 | 9 | 10 |
|---|---|---|---|---|---|---|---|---|---|---|
| Количество сторон многоугольника | 4 | 3 | 3 | 4 | 3 | 6 | 4 | 4 | 4 | 3 |
| Количество итераций до первого возвращения | 2 | 3 | 11 | 20 | 35 | 37 | 63 | 185 | 269 | 479 |

Но заметим, что наряду с большой и маленькой ракетами можно выделить и «среднюю» $А_{ср}$, опирающуюся на стабильную зону номер 3 (рис. 40). Для неё преобразование первого возвращения выглядит гораздо более структурированным, нежели для «маленькой ракеты»; многоугольников оказывается восемь; большие четыре из них возникают естественным образом из разбиения большой ракеты на зоны; поведение же остальных многоугольников очень похоже на поведение больших четырех зон, с той лишь разницей, что маленькие многоугольники перемещаются (рис. 41). Так или иначе, самоподобия с другими рассмотренными ракетами не наблюдается.

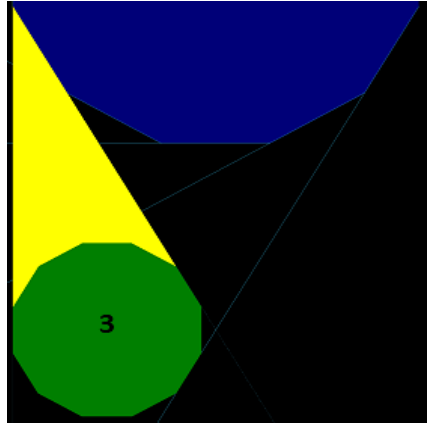

Рис.40. «Средняя ракета» $A_{ср}$

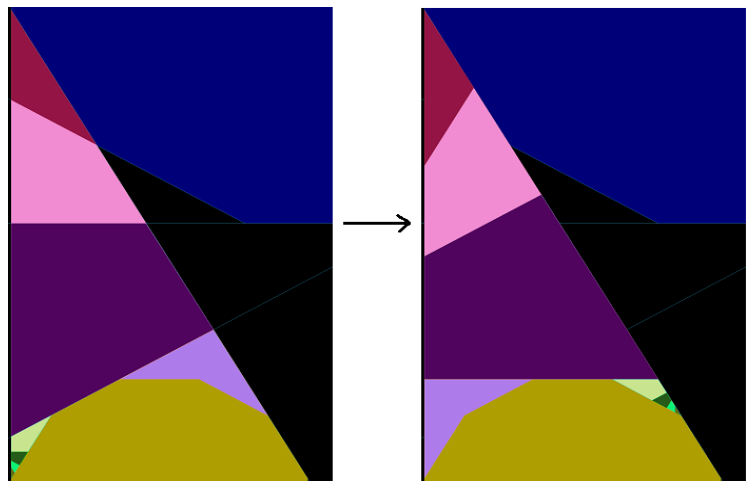

Рис. 41. Преобразование первого возвращения для $A_{ср}$

Однако заметим, что «большую и маленькую ракеты» в механизме можно обобщить до двух подобных фигур, обладающих похожими преобразованиями первого возвращения T. В связи с этим целесообразно рассмотреть «маленькую среднюю» ракету $A_{мал.ср.}$ — ракету, подобную с коэффициентом $\lambda$ средней ракете и имеющую с последней общую вершину. Оказывается, что преобразование первого возвращения T для такой новой ракеты полностью (с точностью до подобия) совпадает с аналогичным преобразованием для средней ракеты! Более формально, если $T_{ср}$ - преобразование первого возвращения для $A_{ср}$, $T_{мал.ср}$ - его аналог для $A_{мал.ср.}$, а Г — сжатие с центром в вершине «ракет» и коэффициентом λ, то для каждой точки $x$ средней ракеты $A_{ср}$, для которых определено преобразование $T_{ср}$, выполняется равенство $T_{мал.ср} Гx = ГT_{ср}x$ !

Перейдем ко второму шагу алгоритма. Чтобы понять, какие точки «средней ракеты» обладают траекториями ранга ноль, т. е. не попадающими в «маленькую среднюю», изобразим траектории точек последней (рис.42):

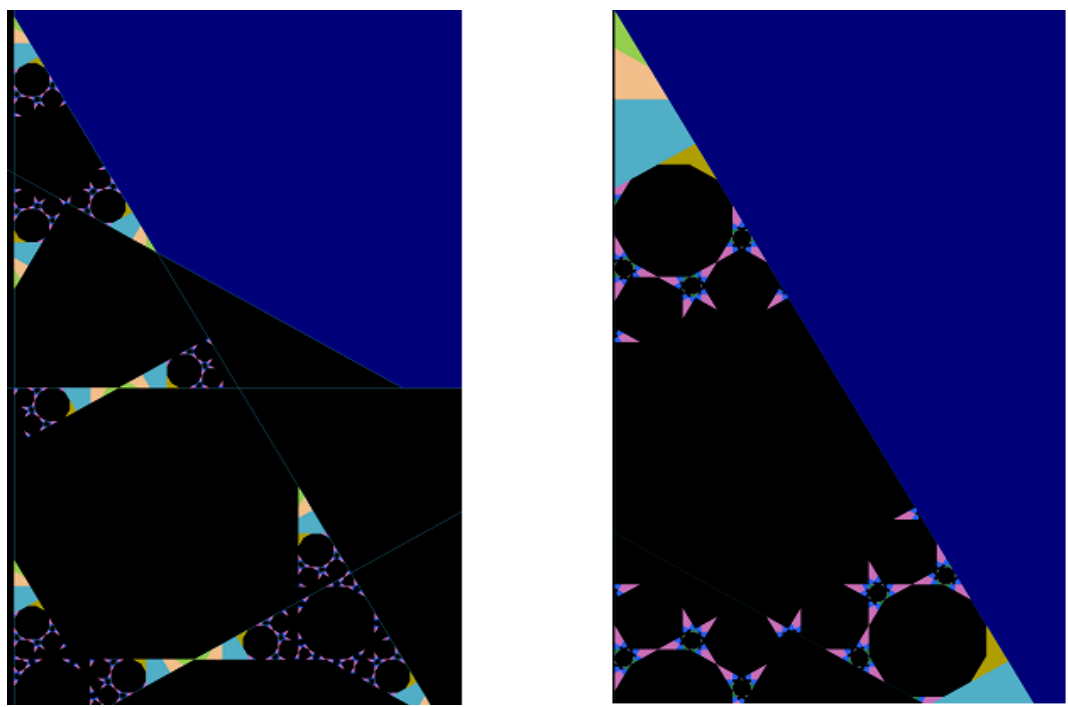

Рис. 42. Траектории точек «маленькой средней» ракеты (справа – в большем масштабе)

Из рисунка можно увидеть некоторые периодические компоненты ранга 0, аналогичные инвариантным зонам; однако для того, что алгоритм сработал, требуется, чтобы вся «средняя ракета» $A_{ср}$ состояла только из выделенных компонент нулевого уровня и следов траекторий точек первого и выше уровней; пустые (черные) зоны на рис.42 не сводятся (по крайней мере, явным образом) к набору выпуклых многоугольников, являющихся периодическими компонентами; это обстоятельство делает невозможным применение второго шага алгоритма.

Однако тот же рисунок 42 дает возможность предположить, что аналогичное самоподобие возникнет и между ракетами A и X. Более точно, пусть $Г_X$ – преобразование, переводящее $A_{ср}$ в X, а $T_X$ есть преобразование первого возвращения T для ракеты X.

**Гипотеза 1:** для всех точек $p$ средней ракеты $A_{ср}$, для которых преобразование $T_{ср}$ определено, выполнено равенство $Г_X T_{ср} p = T_X Г_X p$.

Компьютерные исследования показывают, что гипотеза 1 верна, и количества итераций, требуемые для возвращения, полностью совпадают с количествами итераций, приведенными в таблице 7 приложения. Заметим, что если точка p ракеты $A_{ср}$ периодична относительно T, то p периодична и относительно $T_{ср}$; замечание верно также и для конечных, и для апериодических точек. Аналогично случаю для восьмиугольника можно показать, что $Г_X$ сохраняет тип точек относительно $T_{ср}$.

Докажем теперь существование апериодической траектории. Рассмотрим последовательность периодических многоугольников $C_n$, $n \geq 0$, т.ч. $C_0$, $C_1$, $C_2$ – (открытые) многоугольники, определенные ранее (см. рис.38), а для $n \geq 3$ $C_n = Г_X(C_{n-3})$. Так как $Г_X$ сохраняет периодичность точек, то все $C_n$ суть периодические многоугольники. Согласно гипотезе 1 для любой точки $p$ ракеты $A_{ср}$ выполнено $Г_X T_{ср} p = T_X Г_X p = T_{ср}^{k(p)} Г_X p$, причем $k(p) \geq 2$ (в силу того, что ни одна из точек X при применении $T_{ср}$ не возвращается в X). Отсюда следует, что для всех неотрицательных целых n верно: $per(C_{n+3}) \geq 2*per(C_n)$, где $per(C)$ есть период фигуры C относительно преобразования $T_{ср}$. Таким образом, периоды многоугольников последовательности $C_n$ стремятся к бесконечности; повторяя рассуждения, приведенные в случае правильного восьмиугольника, получаем, что для внешнего биллиарда относительно правильного восьмиугольника существует точка с бесконечной апериодической траекторией. Отметим также, что малая из ракет лежит строго внутри большой, и именно за счет этого факта мы смогли «обойти» пункты алгоритма о полном описании периодических структур ранга 0; однако пункт про самоподобные структуры остается, по-видимому, необходимой и наиважнейшей частью алгоритма.

В заключение отметим, что не только приведенные пары «ракет» могли бы образовывать искомое самоподобие. В приложении можно найти все многоугольники, для которых было проведено исследование преобразование первого возвращения; некоторые пары из них также образуют «хорошее» самоподобие.

## Выводы

Подведем итоги. В данной работе:

1) рассмотрены простые случаи внешнего биллиарда вне правильных решеточных многоугольников; доказано, что в этих случаях не существует бесконечной апериодической траектории;

2) детально исследован внешний биллиард вокруг правильного восьмиугольника; обнаружена ренормализационная схема, исследованы периодические орбиты; построена последовательность точек, пределом которой является точка с бесконечной апериодической траекторией;

3) детально исследован внешний биллиард вокруг правильного двенадцатиугольника; обнаружена (неполная) ренормализационная схема; построена последовательность точек, пределом которой является точка с бесконечной апериодической траекторией.

**Литература**

## Приложение

В данном приложении находятся данные о преобразованиях первого возвращения для различных фигур.

Начнем мы с «самолётика», находящегося над инвариантной фигурой номер один, и его «маленького» аналога, т.е. образа исходного самолётика при применении $Γ_λ$; их преобразования первого возвращения совпадают (с точностью до $Γ_λ$).

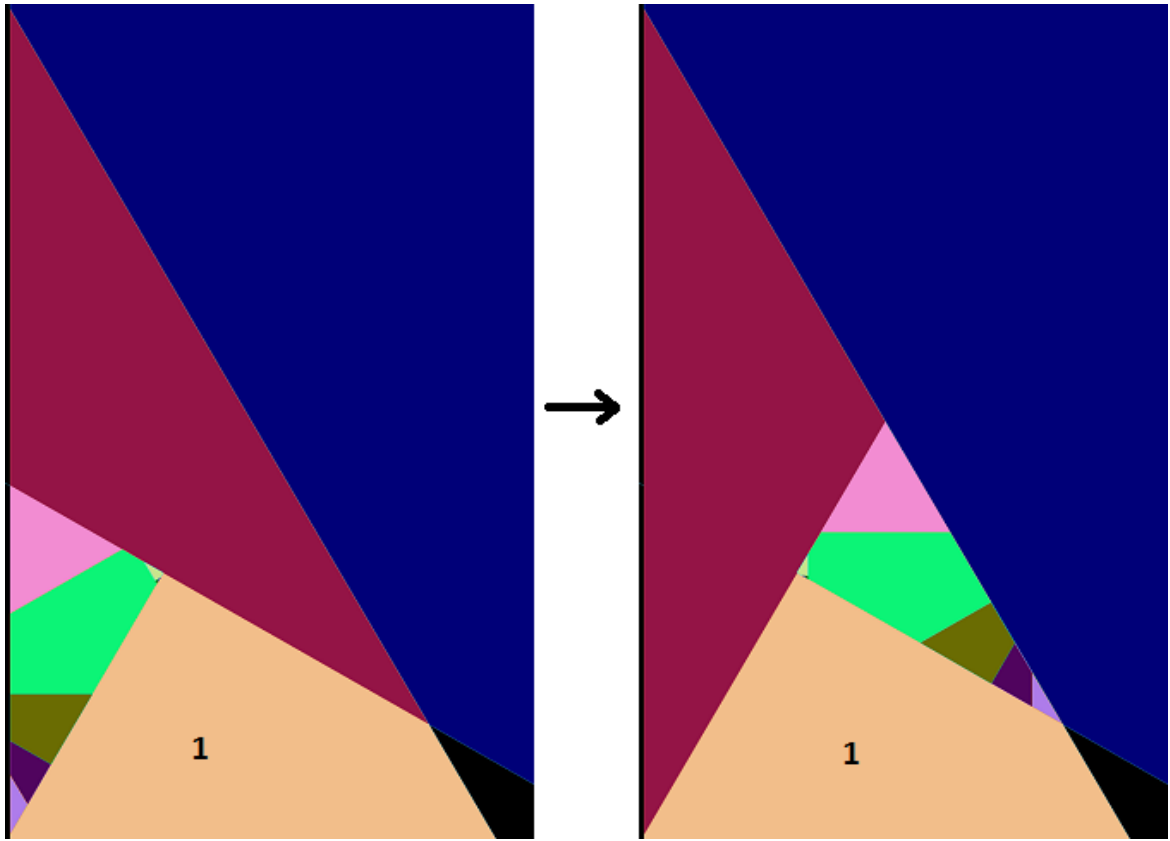

Рис.1. Преобразование первого возвращения для «самолетика»

Таблица 1. Преобразование первого возвращения для «большого самолетика»

|  | 1 | 2 | 3 | 4 | 5 | 6 | 7 | 8 | 9 |
|---|---|---|---|---|---|---|---|---|---|
| Количество сторон многоугольника | 3 | 3 | 4 | 3 | 6 | 4 | 4 | 4 | 3 |
| Количество итераций до первого возвращения | 1 | 9 | 18 | 33 | 35 | 61 | 183 | 267 | 477 |

Таблица 2. Преобразование первого возвращения для «маленького самолетика»

| | 1 | 2 | 3 | 4 | 5 | 6 | 7 | 8 | 9 |
|---|---|---|---|---|---|---|---|---|---|
| Количество сторон многоугольника | 3 | 3 | 4 | 6 | 3 | 4 | 4 | 4 | 3 |
| Количество итераций до первого возвращения | 35 | 429 | 684 | 891 | 1109 | 1353 | 4939 | 7463 | 13773 |

Следующим мы рассмотрим «зону 0» и её «средний» аналог – «среднюю зону 0» (образ зоны 0 при преобразовании большой ракеты в среднюю).

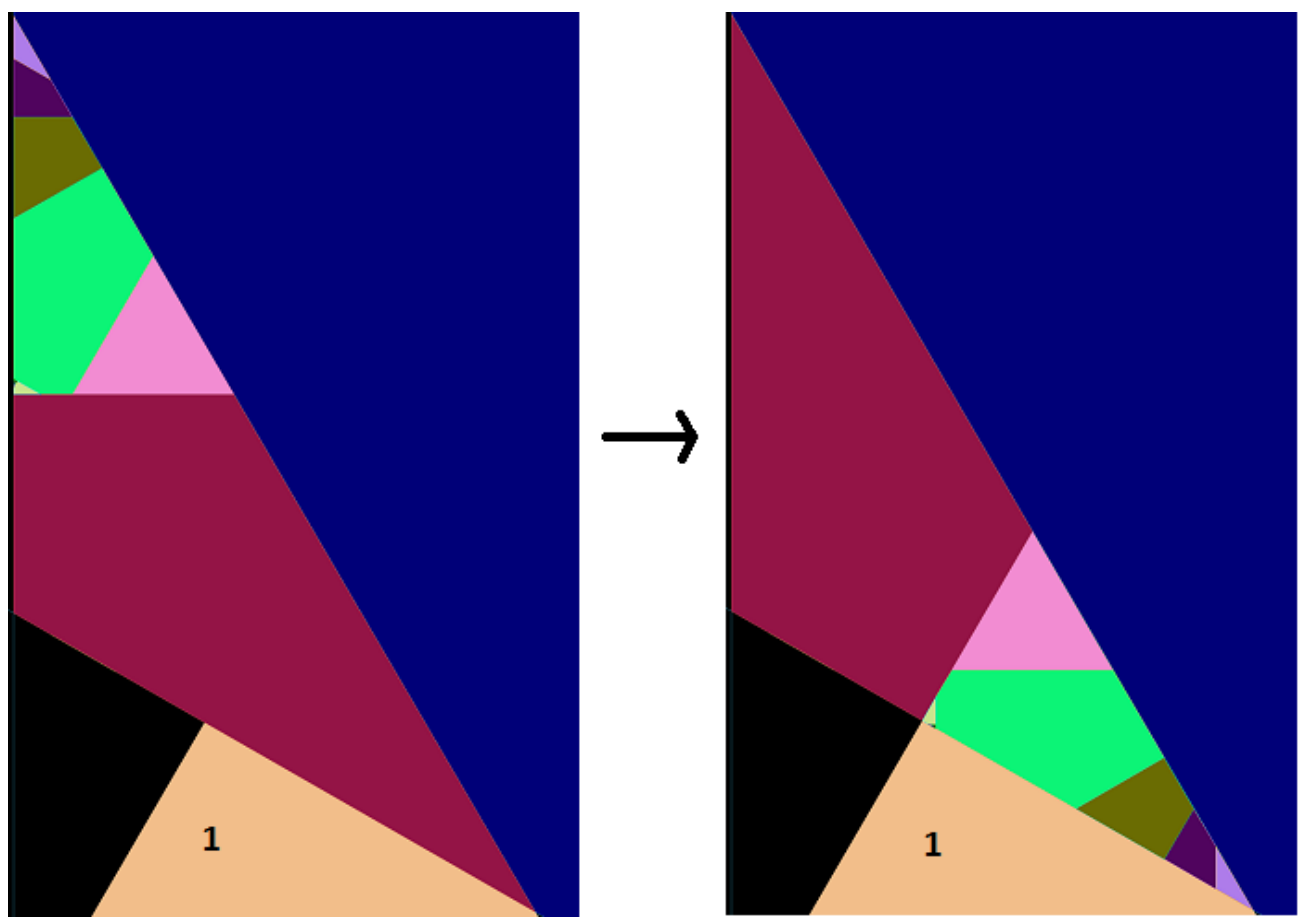

Рис.2. Преобразование первого возвращения для «зоны 0»

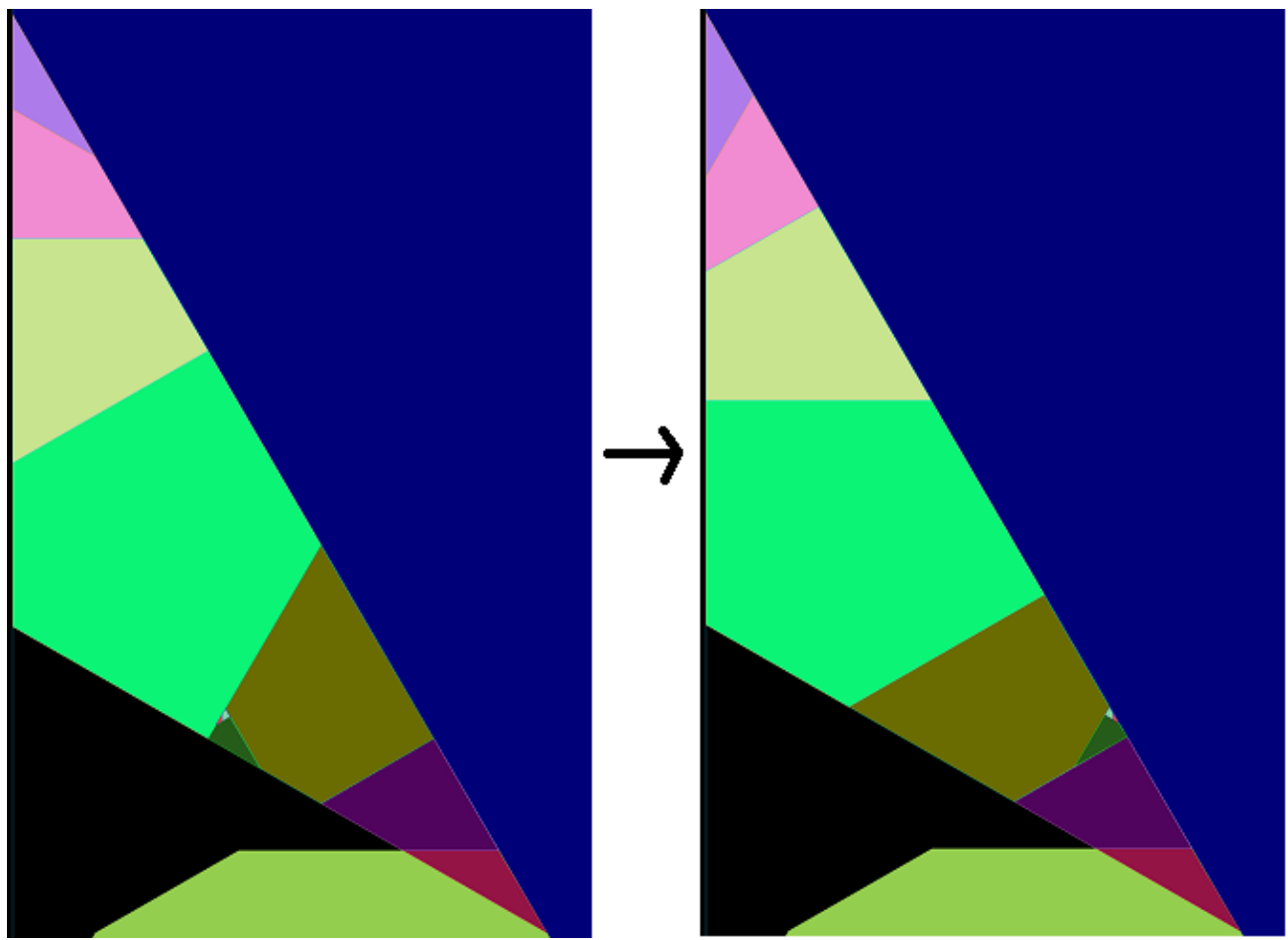

Рис.3. Преобразование первого возвращения для средней «зоны 0»

Таблица 3. Преобразование первого возвращения для «зоны 0»

|  | 1 | 2 | 3 | 4 | 5 | 6 | 7 | 8 | 9 |
|---|---|---|---|---|---|---|---|---|---|
| Количество сторон многоугольника | 4 | 3 | 4 | 3 | 6 | 4 | 4 | 4 | 3 |
| Количество итераций до первого возвращения | 1 | 10 | 19 | 34 | 36 | 62 | 184 | 268 | 478 |

Таблица 4. Преобразование первого возвращения для средней «зоны 0»

|  | 1 | 2 | 3 | 4 | 5 | 6 | 7 | 8 | 9 | 10 |
|---|---|---|---|---|---|---|---|---|---|---|
| Количество сторон многоугольника | 3 | 4 | 4 | 3 | 5 | 5 | 4 | 4 | 4 | 3 |
| Количество итераций до первого возвращения | 18 | 20 | 24 | 35 | 37 | 48 | 63 | 196 | 280 | 490 |

В конце, приведем данные о преобразованиях первого возвращения для различных «ракет». Данные для «маленькой» ракеты приведены ранее; здесь же мы рассмотрим «маленькую-маленькую» ракету (полученную из «маленькой» с помощью $Г_\lambda$); преобразование первого возвращения которой подобно «маленькой» ракете; также в поле нашего зрения попадут преобразования первого возвращения для «средней» и «маленькой средней» (аналогично) ракетам.

Таблица 5. Преобразование первого возвращения для средней «зоны 0»

|  | 1 | 2 | 3 | 4 | 5 | 6 | 7 | 8 | 9 | 10 |
|---|---|---|---|---|---|---|---|---|---|---|
| Количество сторон многоугольника | 4 | 3 | 3 | 4 | 6 | 3 | 4 | 4 | 4 | 3 |
| Количество итераций до первого возвращения | 70 | 105 | 499 | 754 | 961 | 1179 | 1423 | 5009 | 7533 | 13843 |

Таблица 6. Преобразование первого возвращения для средней ракеты $A_{ср}$

|  | 1 | 2 | 3 | 4 | 5 | 6 | 7 | 8 |
|---|---|---|---|---|---|---|---|---|
| Количество сторон многоугольника | 3 | 4 | 4 | 4 | 4 | 3 | 4 | 4 |
| Количество итераций до первого возвращения | 1 | 1 | 1 | 1 | 10 | 25 | 27 | 53 |

Таблица 7. Преобразование первого возвращения для средней ракеты $A_{мал.ср.}$

|  | 1 | 2 | 3 | 4 | 5 | 6 | 7 | 8 |
|---|---|---|---|---|---|---|---|---|
| Количество сторон многоугольника | 4 | 3 | 4 | 4 | 4 | 4 | 3 | 4 |
| Количество итераций до первого возвращения | 20 | 35 | 37 | 63 | 318 | 525 | 743 | 987 |